\documentclass{birkjour}
\usepackage{amssymb}
\usepackage{amsthm}

\usepackage[colorlinks=true, linkcolor=blue, citecolor=blue]{hyperref}
\usepackage{cleveref}
\usepackage{color,graphicx,shortvrb}

\usepackage[active]{srcltx} %SRC Specials for DVI search

\usepackage{enumerate}
\usepackage{amssymb}
\usepackage{tikz-cd}

\usepackage{mathtools}

 \newtheorem{theorem}{Theorem}[section]
 \newtheorem{corollary}{Corollary}[section]
 \newtheorem{lemma}{Lemma}[section]
 \newtheorem{proposition}{Proposition}[section]
 \theoremstyle{definition}
 
 \theoremstyle{remark}
 \newtheorem{remark}{Remark}[section]
 \newtheorem{example}{Example}
 
 \numberwithin{equation}{section}
 \usepackage{amsmath, amssymb, amsthm}
 \usepackage{mathrsfs}
 
\begin{document}

%-------------------------------------------------------------------------
% editorial commands: to be inserted by the editorial office
%
%\firstpage{1} \volume{228} \Copyrightyear{2004} \DOI{003-0001}
%\seriesextraline{This is the Concrete Title of this Book\br H.E. R and S.T.C. W, Eds.}
%\firstpage{1}

\title[Uniqueness of Hanh-Banach extensions and inner ideals]{Uniqueness of Hahn--Banach extensions and inner ideals in real C$^*$-algebras and real JB$^*$-triples}

\author[L. Li]{Lei Li}
\address{School of Mathematical Sciences and LPMC, Nankai University, 300071 Tianjin, China.}
\email{leilee@nankai.edu.cn}

\author[A.M. Peralta]{Antonio M. Peralta}
\address{Instituto de Matem{\'a}ticas de la Universidad de Granada (IMAG), Departamento de An{\'a}lisis Matem{\'a}tico, Facultad de
	Ciencias, Universidad de Granada, 18071 Granada, Spain.}
\email{aperalta@ugr.es}
%\thanks{The research was supported by balabala.}

\author[S. Su]{Shanshan Su}
\address{%
(Current address) School of Mathematics, East China University of Science and Technology, Shanghai, 200237 China.\\
 Departamento de An{\'a}lisis Matem{\'a}tico, Facultad de Ciencias, Universidad de Granada, 18071 Granada, Spain.}
\email{lat875rina@gmail.com}

\author[J. Zhang]{Jiayin Zhang}
\address{%
(Current address) School of Mathematical Sciences and LPMC, Nankai University, 300071 Tianjin, China. \\
Departamento de An{\'a}lisis Matem{\'a}tico, Facultad de Ciencias, Universidad de Granada, 18071 Granada, Spain.}
\email{magicianapo@163.com}
%\thanks{The research was supported by China Scholarship Council.}



\date{today}

\begin{abstract} We show that every closed (resp., weak$^*$-closed) inner ideal $I$ of a real JB$^*$-triple (resp. a real JBW$^*$-triple) $E$ is Hahn--Banach smooth (resp., weak$^*$-Hahn--Banach smooth).  Contrary to what is known for complex JB$^*$-triples, being (weak$^*$-)Hahn--Banach smooth does not characterise (weak$^*$-)closed inner ideals in real JB(W)$^*$-triples. We prove here that a closed (resp., weak$^*$-closed) subtriple of a real JB$^*$-triple (resp., a real JBW$^*$-triple) is Hahn-Banach smooth (resp., weak$^*$-Hahn-Banach smooth) if, and only if, it is a hereditary subtriple. If we assume that $E$ is a reduced and atomic JBW$^*$-triple, every weak$^*$-closed subtriple of $E$ which is also weak$^*$-Hahn-Banach smooth is an inner ideal.\smallskip
	  	
In case that $C$ is the realification of a complex Cartan factor or a non-reduced real Cartan factor, we show that every weak$^*$-closed subtriple of $C$ which is weak$^*$-Hahn-Banach smooth and has rank $\geq 2$ is an inner ideal. The previous conclusions are finally combined to prove the following: Let $I$ 
be a closed subtriple of a real JB$^*$-triple $E$ satisfying the following hypotheses:   

\noindent $(a)$ $I^*$ is separable.
 
\noindent $(b)$ $I$ is weak$^*$-Hahn-Banach smooth. 

\noindent $(c)$ The projection of $I^{**}$ onto each real or complex Cartan factor summand in the atomic part of $E^{**}$ is zero or has rank $\geq 2$. 

\noindent Then $I$ is an inner ideal of $E$. 
\end{abstract}

%%% ----------------------------------------------------------------------
\maketitle
%%% ----------------------------------------------------------------------
%\tableofcontents
\vspace*{-1cm}
\section{Introduction}

Some of the most attractive, and perhaps most surprising, results in mathematics reveal the intrinsic connections between algebraic and analytical properties of different objects. An illustrating example can be the algebraic characterization of $M$-ideals (purely geometric objects) in C$^*$-algebras, JB$^*$-algebras, JB$^*$-triples, and their real forms. More concretely, a projection $P$ on a Banach space $X$ is called an \emph{$M$-projection} (resp., an \emph{$L$-projection}) if  $\|  x \| = max \{ \| P(x) \| , \| (Id-P)(x) \|  \}$ (resp., $\|  x \| = \| P(x) \| + \| (Id-P)(x) \|$), for all $x \in X$. In such a case, we say that $P(X)$ is an \emph{$M$-summand} (resp., \emph{$L$-summand}) of $X$. A closed subspace $M$ of $X$ is said to be an \emph{$M$-ideal} if its polar or annihilator in $X^*$, $M^{\circ}:=\{\varphi\in X^* : \varphi|M \equiv 0\},$ is an $L$-$summand$ of $X^*$ (see \cite{AlfsenEffros1972}). It is only the geometry defined by the norm of $X$ what we need to define $M$-ideals in $X$. There are Banach spaces that admit no $M$-summands but contain an abundant collection of $M$-ideals. The $M$-ideals of a C$^*$-algebra are precisely its closed (two-sided) ideals (see \cite{SmithWard78} or \cite[Theorem V.4.4]{HarmandWernerWernerBook1993}). Subsequent results confirmed the deep interplay between algebra and analysis, for example, the $M$-ideals of the self-adjoint part, $A_{sa}$,  of a C$^*$-algebra $A$ coincide with the subspaces of the form $I\cap  A_{sa},$  with $I$ being a closed two-sided ideal in $A$ \cite[Proposiiton 6.18]{AlfsenEffros1972}, the closed ideals in a JB$^*$-algebra are in one-to-one correspondence with its $M$-ideals \cite{PaPeRo1982}, while in a (complex) JB$^*$-triple, $\mathcal{E}$, the $M$-summands of the Banach space underlying $\mathcal{E}$ are precisely the closed triple ideals in $\mathcal{E}$ \cite[Theorem 3.2]{BartonTimoney1986weak}. In all these results, the arguments rely on the complex linear structure of the involved spaces, so the $M$-ideals in real C$^*$-algebras and real JB$^*$-triples remained undetermined along decades. Quite recently, D. Blecher, M. Neal and the second and third authors of this note finally established the desired characterization in the case of real C$^*$-algebras and real JB$^*$-triples, showing that $M$-ideals of these Banach spaces correspond to (norm) closed ideals and closed triple ideals, respectively \cite{BNPS2025MidealsinrealJB}. Recall that a subspace $\mathcal{I}$ of a (real or complex) JB$^*$-triple $\mathcal{E}$ is called a \emph{triple ideal} if $\{\mathcal{I},\mathcal{E}, \mathcal{E}\}+\{\mathcal{E}, \mathcal{I}, \mathcal{E}\}\subseteq \mathcal{I}$.\smallskip

On the other hand, an interesting property of $M$-ideals assures that every bounded functional in the dual space, $M^*$, of an $M$-ideal $M$ in a Banach space $X$ has a unique norm-preserving extension to a functional in $X^*$ (cf. \cite[Proposition I.1.12]{HarmandWernerWernerBook1993}). So, it is natural to ask what are the closed subspaces satisfying this unique Hahn--Banach extension property. Following \cite[page 44]{HarmandWernerWernerBook1993} (see also \cite{SmithSullivan1976,Sullivan1977}), we shall say that a closed subspace $Z$ of a Banach space $X$ is \emph{Hahn--Banach smooth} if every functional in $Z^{*}$ admits a unique norm-preserving extension to a functional in $X^*$. If $Z$ is a weak$^*$-closed subspace of a dual Banach space $X$, with predual $X_*$, we shall say that $Z$ is \emph{weak$^*$-Hahn--Banach smooth} if every functional in $Z_{*}$ admits a unique Hahn--Banach norm-preserving extension to a functional in $X_*$. As we shall see next, in some well known structures, including C$^*$-algebras and (complex) JB$^*$-triples, the subalgebras and subtriples which are Hahn--Banach smooth can be algebraically characterised.\smallskip 

Let us call to mind that an inner ideal in a real or complex JB$^*$-triple $\mathcal{E}$ is a closed subspace $\mathcal{I}$ satisfying $\{\mathcal{I},\mathcal{E},\mathcal{I}\}\subseteq \mathcal{I}$ (note that along this note all inner ideals will be assumed to be closed). For example, in a C$^*$-algebra $A$, subspaces of the form $p A q$ are inner ideals of $A$, which are not, in general, (triple) ideals. A fascinating result by C.M. Edwards and G.T. Rüttimann asserts that a norm-closed subtriple $\mathcal{B}$ of a complex JB$^*$-triple $\mathcal{E}$ is an inner ideal in $\mathcal{E}$ if, and only if, every bounded linear functional on $\mathcal{B}$ has a unique norm-preserving linear extension to $\mathcal{E}$ \cite{EdwRutti1992}. Despite  Edwards and Ruttimann addressed the study of inner ideals in real JBW$^*$-triples in \cite{EdwRutti2001}, the problem whether norm closed inner ideals in C$^*$-algebras and real JB$^*$-triples can be characterised as those JB$^*$-subtriples for which every continuous functional admits a unique Hahn--Banach extension remained as an open problem  (see section~\ref{sec: preliminaries} for the proper definitions). As observed by Goodearl in his book \cite{Goodearl1982},  \emph{``The change of coefficient field from $\mathbb{C}$ to $\mathbb{R}$ is more than just a cosmetic change.''} The just posed question is an illustrative example, and has remained open for years.\smallskip 

In this note we shed some new light on our knowledge on Hahn-Banach smooth closed subtriples of real JB$^*$-triples. In a first result we show that every inner ideal in a real JB$^*$-triple is Hahn-Banach smooth (cf. Corollary~\ref{c characterization of closed inner ideals}), and every weak$^*$-closed inner ideal in a real JBW$^*$-triple is weak$^*$-Hahn-Banach smooth (see Theorem~\ref{t characterization of wstarHahnBanach smooth subtriples}). The divergences with respect to the known results for complex JB$^*$-triples are not long in appearing. As we shall see in Example~\ref{example w*HahnBanach smoothness is not enough}, there exists finite dimensional subtriples of real JBW$^*$-triples which are (weak$^*$-)Hahn-Banach smooth but fail to be inner ideals. Besides the characterization of all  weak$^*$-Hahn-Banach smooth weak$^*$-closed subtriples of real JBW$^*$-triples obtained in the mentioned Theorem~\ref{t characterization of wstarHahnBanach smooth subtriples}, a more elaborated argument allows us to characterize all Hahn-Banach smooth closed subtriples $I$ of real JB$^*$-triple $E$ in terms of Peirce subspaces of range tripotents of elements in $I$ in the second dual space of $E$. Actually, the combination of theorems \ref{t characterization of closed subtriples with HB smoothness} and \ref{t characterization of closed subtriples with weak*HB smoothness as hereditary subtriples} leads to the following: Let $I$ be a closed subtriple of a real JB$^*$-triple $E$. Then the following statements are equivalent:\begin{enumerate}[$(a)$]\item $I$ is Hahn--Banach smooth.
	\item $\displaystyle I =\bigcup_{a\in I, \|a\| =1} (E^{**})^{1}(r(a))\cap E,$ where $r(a)$ denotes the range tripotent of $a$ in $E^{**}$ and $(E^{**})^{1}(r(a))$ is the Peirce subspace of all $x\in E^{**}$ satisfying $\{r(a),x,r(a)\} = x$.
	\item $I$ is a hereditary subtriple of $E$, that is, for each $a\in F$, the inner ideal of $I$ generated by $a$, $I(a)$, is a hereditary JB$^*$-subalgebra of of the inner ideal, $E(a)$, of $E$ generated by $a$, equivalently, if $b\in F(a)$ and $c\in E(a)$ are two positive elements with $c\leq b$ in $E(a)$, we have $c\in F$.
\end{enumerate}

Section~\ref{sec: reduced Cartan factors} is devoted to studying those weak$^*$-closed subtriples of reduced atomic real JBW$^*$-triples which are weak$^*$-Hahn-Banach smooth. We recall that a real Cartan factor $C$ is called reduced if $C^{-1} (e) =\{0\}$ for every minimal tripotent $e\in C$, and a real JBW$^*$-triple is called reduced and atomic if it coincides with a direct sum of reduced real Cartan factors. The conclusion in this case is closer to what we know for complex JBW$^*$-triples. Let $W$ be a reduced and atomic real JBW$^*$-triple, and suppose that $I$ is a weak$^*$-closed subtriple of $W$. We establish in Theorem~\ref{thm reduced  atomic JBW-star triples} the equivalence of the following statements:
\begin{enumerate}[$(a)$]
	\item $I$ is an inner ideal.
	\item $I$ is weak$^*$-Hahn--Banach smooth.
\end{enumerate}

The case of complex Cartan factors deserves its own attention in section~\ref{sec: complex Cartan factors}. The counterexamples exhibited in section~\ref{sec: inner ideals and hereditary subtriples} point out that we need to consider extra hypotheses in this setting. Our main conclusion, established in Theorem~\ref{t real subtriples of complex Cartan factors}, asserts the following: Let $I$ be a weak$^*$-closed real subtriple of a complex Cartan factor $C$. Suppose $I$ satisfies the following hypotheses:\begin{enumerate}[$(a)$]
	\item $I$ is weak$^*$-Hahn-Banach smooth.
	\item $I$ has rank bigger than or equal to $2$. 
\end{enumerate} Then $I$ is a complex subtriple and an inner ideal of $C$.\smallskip

Non-reduced real Cartan factors of rank-one always admit weak$^*$-Hahn-Banach smooth weak$^*$-closed subtriple which are not inner ideals. In case that $C$ is a non-reduced real Cartan factor, and $I$ is a weak$^*$-closed subtriple of $C$ having rank $\geq 2$, we prove that $I$ is weak$^*$-Hahn-Banach smooth if, and only if, $I$ is an inner ideal (cf. Propositions~\ref{p type I2p2q{{H}}} and \ref{p III2p{H}}).\smallskip

The previous conclusions are combined in Theorem~\ref{t sufficient conditions for HBS implies inner ideal} to establish the next result: Let $I$ be a closed subtriples of a real JB$^*$-triple $E$. Suppose that $I$ satisfies the following hypotheses:
\begin{enumerate}[$(a)$]\item $I^*$ is separable.
	\item $I$ is Hahn-Banach smooth.
	\item The projection of $I^{**}$ onto each real or complex Cartan factor summand in the atomic part of $E^{**}$ is zero or has rank greater than or equal to $2$.
\end{enumerate} Then $I$ is an inner ideal of $E$.\smallskip

It is worth recalling at this stage that, contrary to the complex setting, surjective linear isometries between real JB$^*$-triples need not preserve triple products. However, if $T: E\to F$ is a surjective linear isometry between two real JB$^*$-
triples such that $E^{**}$ does not contain (real or complex) rank-one Cartan factors as a summand, then $T$ is a triple isomorphism (cf. \cite[Theorem 3.2]{PoloPe2004surjective}).\smallskip

The paper also contains a post-credits section~\ref{sec: post-credits}, where we present a simplified argument to establish that every $M$-summand (resp., every $M$-ideal) in a real JBW$^*$-triple (resp., in a real JB$^*$-triple) is a triple ideal (see Theorem~\ref{t every M-summand is a weak*-closed subtriple}), which just uses a tool based on the facial structure of the closed unit ball of a real JBW$^*$-triple, and reduces significantly the arguments in \cite{BNPS2025MidealsinrealJB}.   

\subsection{Notation and preliminaries}\label{sec: preliminaries}

C$^*$-algebras are among the best known and studied mathematical objects in mathematics. Their structure, classification and geometric properties have been deeply studied since the decade of 1940's. C$^*$-algebras are included inside the strictly wider class of complex Banach spaces known as JB$^*$-triples. It is perhaps worth recalling that a (complex) \emph{JB$^*$-triple} $\mathcal{E}$ is a complex Banach space endowed with a continuous triple product $\{ \cdot, \cdot, \cdot \}$ from $\mathcal{E} \times \mathcal{E} \times \mathcal{E}$ into $\mathcal{E}$, which is symmetric in the first and third position and conjugate linear in the second position, and satisfies:
\begin{enumerate}[$(i)$]
\item For every $x,y,z,a,b \in \mathcal{E}$, we have 
$$\{a,b, \{x,y,z \} \}= \{ \{a,b,x\},y,z \} - \{x, \{b,a,y\},z \} +\{x,y, \{a,b,z \} \} ;$$ \hfill (\emph{Jordan identity})
\item For every $a\in \mathcal{E}$, the linear operator $L(a,a)$ is hermitian with non-negative spectrum;
\item $\| \{a,a,a \} \| = \|a\| ^{3}$ for all $a \in \mathcal{E}$. \hfill (\emph{Gelfand-Naimark axiom})	
\end{enumerate}
Given $x,y \in \mathcal{E}$, we write $L(x,y)$ and $Q(x,y)$ for the linear and conjugate linear maps on $\mathcal{E}$ defined by $L(x,y) (z) = \{ x,y,z \}$ and $Q(x,y) (z) = \{ x,z,y \}$ for all $z \in \mathcal{E}$, respectively.\smallskip

A (complex) \emph{JBW$^*$-triple} is a JB$^*$-triple which is also a dual Banach space. There are some basic properties of JBW$^*$-triples. For instance, the triple product is separately weak$^*$-continuous in every JBW$^*$-triple; the bidual of each JB$^*$-triple is a JBW$^*$-triple and each JBW$^*$-triple admits a unique isometric predual. The reader is referred to \cite{FriedmanRusso1985structure, Upmeier1985jordan, BartonTimoney1986weak,Dineen1986second,EdwRutti1992} and the references therein for the basic theory of (complex) JB$^*$-triples and JBW$^*$-triples.\smallskip

In this paper we shall be mainly interested in real C$^*$-algebras and real JBW$^*$-triples. We recall that a real C$^*$-algebra ${A}$ is a closed $^*$-invariant real subalgebra of a (complex) C$^*$-algebra (see \cite[Definition 5.1.1 and Proposition 5.1.2]{BingrenBook2003}). The symbol ${A}_{sa}$ will stand for the real subspace of all hermitian or self-adjoint elements in $A$. It is known that $A_{sa}$ is not, in general, a subalgebra of $A$ for the associative product. However $A_{sa}$ is a JB-algebra in the usual sense (see \cite{HanchOlsenStormerBook}) when we consider the natural Jordan product defined by $a\circ b := \frac12 (a b + ba )$. We shall write ${A}^+$ for the set of all positive elements in ${A}$. \smallskip

According to the foundational reference \cite{IsidroKaupPalacios1995realform}, a \emph{real JB$^*$-triple} is a closed real subtriple of a complex JB$^*$-triple. Real JB$^*$-triples can be equivalently defined as ``real forms'' of complex JB$^*$-triples, concretely, $E$ is a real JB$^*$-triple if there exists a complex JB$^*$-triple $\mathcal{E}$ and a conjugation (i.e. a conjugate linear isometry of period 2) $\tau$ on $\mathcal{E}$ such that 
$$E = \mathcal{E}^{\tau} = \{ x \in \mathcal{E}: \tau(x) = x \}.$$  The characterization is better understood if we have in mind that every surjective linear or conjugate-linear isometry on a JB$^*$-triple is a triple isomorphism (cf. \cite[Proposition 5.5]{Kaup1983riemann} and \cite[Definition 2.1 and subsequent comments]{IsidroKaupPalacios1995realform}). If $ \mathcal{E}$ is a C$^*$-algebra and $\tau$ is a conjugate-linear $^*$-automorphism or period 2 on $E$, the real form $\mathcal{E}^{\tau}$ is a real C$^*$-algebra and all real C$^*$-algebras are obtained in this way (see \cite[Proposition 5.1.3]{BingrenBook2003} and \cite{Goodearl1982}).\smallskip

As in the case of complex JBW$^*$-triples, a \emph{real JBW$^*$-triple} is a real JB$^*$-triple which is also a dual Banach space. Analogously to the complex case, every real JBW$^*$-triple has a unique isometric predual and its triple product is separately weak$^*$-continuous \cite{MartinezPe2000separate}. The bidual, $E^{**},$ of a real JB$^*$-triple, $E$, is a real JBW$^*$-triple \cite[\S 4]{IsidroKaupPalacios1995realform}.
\smallskip

It should be noted that the class of real JB$^*$-triples is extremely wide, it strictly contains the classes of real and complex C$^*$-algebras, complex JB$^*$-triples, JB-algebras, real JB$^*$-algebras, real and complex spin factors and Hilbert spaces, and real forms of exceptional Cartan factors (cf. \cite{IsidroKaupPalacios1995realform,Kaup1997,MartinezPe2000separate, EdwRutti2001}).\smallskip

We recall that a real JB$^*$-algebra $\mathfrak{A}$ is a closed $^*$-invariant real subalgebra of a (complex) JB$^*$-algebra in the sense of Kaplansky \cite{Wright1977, HanchOlsenStormerBook}. We shall write $\mathfrak{A}_{sa}$ for the JB-algebra of all self-adjoint elements in $\mathfrak{A}$, while $\mathfrak{A}_{skew}$ will denote the set of all skew-symmetric elements in $\mathfrak{A}$. Unital JB$^*$-algebras were introduced by K. Alvermann under the name J$^*$B-algebras (see \cite{Alvermann1986}). It follows from the definition that every JB-algebra is a real JB$^*$-algebra. An element $a$ in $\mathfrak{A}$ is called positive if it is self-adjoint and its spectrum is contained in $\mathbb{R}^{+}$. The set of all positive elements in $\mathfrak{A}$ will be denoted by $\mathfrak{A}^+$. \smallskip

Let $E$ be a (real or complex) JB$^*$-triple. A closed subtriple $I \subseteq E$ is called an \textit{inner ideal} of $E$ if $\{I,E,I \} \subseteq I,$ and it is called a \emph{(triple) ideal} of $E$ if $\{E,E,I\} + \{E,I,E \} \subseteq I$. Inner ideals and ideals have been studied from an algebraic point of view in the case of real and complex C$^*$-algebras, and real and complex JB$^*$-triples (see \cite{BartonTimoney1986weak, Horn1987characterization, EdwRutt1989inneridearinWalgebras,EdwRutt1991inneridearinCalgebras,EdwRutti2001,BNPS2025MidealsinrealJB, BNPSMidealsOperatorAlgebras} for more details). \smallskip

An element $e$ in a real or complex JB$^*$-triple $E$ is called a \emph{tripotent} if $\{e,e,e\} = e$. The symbol $\mathcal{U}(E)$ will stand for the set of all tripotents in $E$. Each tripotent $e \in E$ induces two decompositions of $E$ as follows:
$$ E = E_{0} (e) \oplus E_{1} (e) \oplus
E_{2} (e) = E^{0} (e) \oplus  E^{1} (e) \oplus E^{-1} (e) $$
where $$E_{k} (e) := \Big\{ x\in E : L(e,e)x = \frac{k}{2} x \Big\} \ \ \ (k = 0,1,2)$$ which is a subtriple of $E,$ and $$E^{k} (e) := \left\{ x\in E : Q(e) (x) := \{ e,x,e \}= k x \right\} \ \ \ (k = 0,1,-1)$$ is a Banach subspace of $E$. The first decomposition is known as the \emph{Peirce decomposition} of $E$ associated with $e.$ Henceforth, the natural projection of $E$ onto $E_{k} (e)$ (resp., onto $E^{k} (e)$) will be denoted by $P_{k} (e)$ (resp., $P^{k} (e)$). The projection $P_k(e)$ is known as the \emph{Peirce-$k$ projection}. It is known that all Peirce projections associated with a tripotent $e$, as well as $P^{1} (e),$ are contractive (see \cite[Remark 2.6]{PeStacho2001}). The previous decompositions of $E$ obey the following multiplication arithmetic:
\begin{equation*}
\textit{Peirce rules}\left\{
\begin{aligned}
& \{ {E_{i}(e)},{E_{j}(e)},{E_{k}(e)} \} \subseteq E_{i-j+k}(e),\ \hbox{
if  } i,j,k\in\{0,1,2\}, \\ 
& \{ {E_{i}(e)},{E_{j}(e)}, {E_{k}(e)} \} =\{0\}, \hbox{ for } i-j+k\neq 0,1,2, \\
& \{ {E_0 (e)},{E_2 (e)},{E}\} = \{ {E_2 (e)},{E_0 (e)},{E }\} =\{0\};
\end{aligned}
\right.
\end{equation*}
and
$$\begin{aligned}
		& E_{2} (e) = E^{1} (e) \oplus E^{-1} (e), \quad E_{1} (e) \oplus E_{0} (e) = E^{0} (e), \\
	   	& \{ {E^{i} (e)},{E^{j} (e)},{E^{k} (e)}\} \subseteq E^{i j k} (e),
		\hbox{ whenever } i j k \ne 0. \ \ \ \ \ \ \ \ \ \ \ \nonumber
\end{aligned}$$
Furthermore, when equipped with the Jordan product and involution defined by $x\circ_{e} y := \{x,e,y\}$ and $x^{{*}_{e}} := \{e,x,e\}$, respectively, the Peirce-$2$ subspace $E_{2}(e)$ becomes a real or complex unital JB$^*$-algebra with identity $e$. The self-adjoint part of this real or complex JB$^*$-algebra coincides with  $E^{1}(e),$ and hence the latter is a unital JB-algebra. The Peirce subspaces $E_2(e)$ and $E_0(e)$ are inner ideals of $E$. One of the challenges that might occur when we investigate a problem on a general real JB$^*$-triple $E$ is that, for a tripotent $e\in E$, the JB$^*$-subtriples $E^{-1}(e)$ and $E^{1}(e)$ are not, in general, related (while in a complex JB$^*$-triple $\mathcal{E}$ we have $\mathcal{E}^{-1}(e) = i \mathcal{E}^{1}(e)$).\smallskip

The projections $P_k(e)$ and $P^{k}(e)$ can be described in terms of the triple product via the following identities:$$\begin{aligned}\label{expresions for Peirce projections}
P_2 (e) &= Q(e)^2, \  P_1 (e) = 2 L(e,e) -2 Q(e)^2, \  P_0(e) = Id_{E} - 2 L(e,e) + Q(e)^2, \\
P^{1} (e) &= \frac12 \left(Q(e)^2 + Q(e)\right), \ P^{-1} (e) = \frac12 \left(Q(e)^2 - Q(e)\right), \hbox{ and } \\  P^{0} (e) &= P_0(e)+P_1(e).
\end{aligned}$$ In case that $e$ is a tripotent in a real JBW$^*$-triple $W$, it follows from the separate weak$^*$-continuity of the triple product of $W$ that $P^{j} (e)$ and $P_k(e)$ are weak$^*$-continuous\label{ref Peirce projections and upper projections are weakstar conts} projections for all $k\in\{0,1,2\}$ and $j\in \{-1,1,0\}$.\smallskip

A non-zero tripotent $e$ in a real or complex JB$^*$-triple $E$ is called \emph{minimal} if $E^{1} (e) = \mathbb{R} e$. In case that $E$ is a complex JB$^*$-triple, $e$ is minimal if, and only if, $E_2 (e) = \mathbb{C} e$.\smallskip

If $A$ is a real C$^*$-algebra regarded as a real JB$^*$-triple with product $\{a,b,c\}= \frac12 (a b^* c+ cb^*  a),$ tripotents in $A$ coincide with partial isometries. The Peirce decompositon of $A$ associated with a partial isometry $e$ is given by $A_2(e) = ee^* A e^*e,$ $A_1(e) = (\mathbf{1} -ee^*) A e^*e \oplus ee^* A (\mathbf{1}-e^*e)$, and  $A_0(e) = (\mathbf{1} -ee^*) A (\mathbf{1}-e^*e)$.\smallskip

Elements $x,y$ in a real JB$^*$-triple $E$ are said to be \emph{orthogonal} ($x \perp y$ in short) if, and only if, $L(x,y) = 0$ if and only $\{x,x,y\} = 0$ (see \cite[Lemma 1]{BurgosPoloGarcesPe2008orthogonality} for more equivalent characterizations). Furthermore, two JB$^*$-subtriples $M$ and $N$ in $E$ are orthogonal (denoted as $M\perp N$) if $\{m,n,E\} = \{0\}$ for all $m\in M$, $n\in N$, or equivalently, $L(M,N) =0$.\smallskip

We conclude this brief background section by recalling the notion of continuous triple functional calculus. Suppose now that $\mathcal{E}$ is a complex JB$^*$-triple, and $x$ is a fixed element in $\mathcal{E}$. Denote $x^{[1]} := x$, $x^{[3]}:= \{x,x,x\}$ and $x^{[2n+1]} := \{x,x^{[2n-1]}, x\} $ for $n \in \mathbb{N}$.
Let the symbol $\mathcal{E}_{x}$ stand for the JB$^*$-subtriple of $\mathcal{E}$ generated by $x$. It is part of JB$^*$-triple theory that $\mathcal{E}_{x}$ is (triple) isometrically isomorphic to $C_{0}(\Omega_{x})$ for some locally compact Hausdorff space $\Omega_{x} \subseteq (0, \|x\| ]$ with $\Omega_{x} \cup \{0\}$ compact. Moreover, there is a triple isomorphism $\Psi:\mathcal{E}_{x}\to C_{0}(\Omega_{x})$ satisfying $\Psi(x)(t) = t$ for $t \in \Omega_{x}$ (cf. \cite[Corollary 1.15]{Kaup1983riemann} and \cite[Corollary 4.8]{Kaup1977algebraic}). For each continuous function $f\in C_{0}(\Omega_{x})$, the \emph{continuous triple functional calculus of $f$ at $x$} is defined as $f_{t} (x) = \Psi^{-1} (f).$ We can define in this way the $(2n-1)$th-root of the element $x$ ($n\in \mathbb{N}$), which is denoted by $x^{[1/(2n-1)]}$ and satisfies $(x^{[1/(2n-1)]})^{[2n -1]} = x$.\smallskip

According to notation in the previous paragraph, for each norm-one element $x$ in a (complex) JB$^*$-triple $\mathcal{E}$, we can define the sequences $(x^{[2n-1]})_{n}$ and $(x^{[1/(2n-1)]})_{n}$. In case that $\mathcal{E}$ is a JBW$^*$-triple, we know that these sequences converge in the weak$^*$-topology of $\mathcal{E}$ to certain tripotents $u(x)$ and $r(x)$ in $\mathcal{E}$, known as the \emph{support} and the \emph{range tripotent} of $x$ in $\mathcal{E}$, respectively (cf. \cite[\S 3]{EdwardsRutti1988facial}). %The reader should be warned that in \cite{EdwardsRutti1996compact}, $r(x)$ is called the support tripotent of $x$. 
The support and range tripotents satisfy $u(x) \leq x \leq r(x)$ in the local order given by the cone of positive elements in the JBW$^*$-algebra $\mathcal{E}_{2}(r(x))$. The range tripotent of $x$ in $\mathcal{E}$ is actually the smallest tripotent $e$ in $\mathcal{E}$ satisfying that $x$ is positive in $\mathcal{E}_{2}(e)$ (see \cite[Lemma 3.6]{EdwardsRutti1988facial}). If $\mathcal{E}$ is a mere JB$^*$-triple, the range tripotent of each norm-one element $x\in \mathcal{E}$ is computed in $\mathcal{E}^{**}$.

Suppose now that $W$ is a real JBW$^*$-triple. It is known from \cite{MartinezPe2000separate} that there exist a JBW$^*$-triple $\mathcal{W}$ and a weak$^*$-continuous conjugate-linear triple automorphism $\tau$ on $\mathcal{W}$ satisfying $W = \mathcal{W}^{\tau}$. Tripotents in $W$ correspond to $\tau$-symmetric tripotents in $\mathcal{W}$. Furthermore, for each norm-one element $x$ in $W$ the range and support tripotents of $x$ in $\mathcal{W}$ are $\tau$-symmetric, and hence they both belong to $W$, that is, $\tau (u_{_\mathcal{W}}(x)) = u_{_\mathcal{W}}(x)\in W$, $\tau (r_{_\mathcal{W}}(x)) = r_{_\mathcal{W}}(x)\in W$. We write $r_{_W} (x) = r_{_\mathcal{W}} (x)$ and $r_{_W} (x) = r_{_\mathcal{W}} (x)$, and we call them the support and range tripotents of $x$ in $W$, respectively (cf. \cite{EdwRutti2001} and \cite[\S 3]{CuetoPeralta2019}).\smallskip 

Observe that if $V$ is a weak$^*$-closed real JB$^*$-subtriple of a real JBW$^*$-triple $W$, for each norm-one element $x$ in $V$ we have $r_{_V} (x) = r_{_W} (x)$ since both tripotents can be computed as weak$^*$-limits of the sequence $(x^{[1/(2n-1)]})_{n}$ in $V$ and $W$, respectively, and $V$ is weak$^*$-closed. A similar conclusion holds for the support tripotent of $x$.\smallskip

We shall finally focus on some special inner ideals. Given a non-zero element $x$ in a complex JB$^*$-triple $\mathcal{E}$, we shall write $\mathcal{E}(x)$ for the inner ideal generated by the element $x$ in $\mathcal{E}$. The structure of $\mathcal{E}(x)$ has been perfectly determined in \cite[Proposition 2.1]{Bunce2000structure}. It follows from the just quoted result that $\mathcal{E}(x)$ coincides with the norm-closure of $\{x, \mathcal{E}, x\} = Q(x) (\mathcal{E})$ in $\mathcal{E}$, and satisfies the following properties:
\begin{enumerate}[$(1)$] 
	\item $\mathcal{E}_{x} \subset \mathcal{E}(x),$
	\item $\mathcal{E}(x)$ is a JB$^*$-subalgebra of $\mathcal{E}^{**}_2(r(x))$,
	\item  $\mathcal{E}(x)^{**} = \overline{\mathcal{E}(x)}^{w^*} = \mathcal{E}_{2}^{**}(r(x)),$ where we write $\overline{\mathcal{E}(x)}^{w^*}$ for the weak$^*$-closure of $\mathcal{E}(x)$ in $\mathcal{E}^{**}$.
\end{enumerate}
 It is further known that $\mathcal{E}^{**}_2 (r(x)) \cap \mathcal{E}$ is weak$^*$-dense in $\mathcal{E}^{**}_2 (r(x))$ (cf. \cite[page 167]{EdwardsRutti1996compact}). Observe that the JB-algebra $\mathcal{E}(x)_{sa}$ is contained in $(\mathcal{E}^{**})^{1} (r(x)) \cap \mathcal{E}$ and is weak$^*$-dense in $(\mathcal{E}^{**})^{1} (r(x))$. 
It follows from these comments that, for each norm-one element $a$ in a real JB$^*$-triple $E$, passing through the complexification, $\mathcal{E}$, of $E$, the following properties hold: \begin{enumerate}[$(1)$]\label{eq density properties of the self-adjoint part}
\item The space $E(a)= \overline{Q(a)(E)}$ is a real JB$^*$-subalgebra of $\mathcal{E}(a)$.
\item $E(a)$ is weak$^*$-dense in $E^{**}_2 (r(a)).$
\item The JB-algebra $E(a)_{sa}$ is contained in $(E^{**})^1 (r(a))\cap E$ and is weak$^*$-dense in $(E^{**})^1 (r(a)).$
\end{enumerate}

\section{Inner ideals and hereditary subtriples}\label{sec: inner ideals and hereditary subtriples}

Our first proposition is devoted to establish an algebraic characterization of weak$^*$-closed inner ideals in real JBW$^*$-triples that will be crucial for later purposes.

\begin{proposition}\label{algebraic charac of inner ideals}
Let $W$ be a real JBW$^*$-triple, and let $ I \subseteq W $ be a weak$^*$-closed subtriple. Then, the following conditions are equivalent:
\begin{enumerate}[$(i)$]
\item $I$ is an inner ideal in $W$.
\item $\displaystyle I =\bigcup_{e\in \mathcal{U}(I)} W^1(e)  = \bigcup_{e\in \mathcal{U}(I)} W_2(e)$.
\end{enumerate}
\end{proposition}

\begin{proof}$(i)\Rightarrow (ii)$ Suppose $I$ is an inner ideal in $W$. It is clear that $$\displaystyle \bigcup_{e\in \mathcal{U}(I)} W^1(e)\subseteq \bigcup_{e\in \mathcal{U}(I)} W_2(e),$$ since $W^1(e)\subseteq W_2(e)$ for every $e\in \mathcal{U}(I)$. Moreover, by the identities in page~\pageref{expresions for Peirce projections} $W_2 (e) = P_2 (e) (W) = Q(e)^2 (W) = Q(e) Q(e)(W) \subseteq I$, because $I$ is an inner ideal and $e$ lies in $I$. It is easy to see that $I\subseteq \bigcup_{e\in \mathcal{U}(I)} W^1(e)$ because every $a\in I$ lies in $W^1(r_{_I} (a))$.  \smallskip 
        
\noindent $(ii) \Rightarrow (i)$ Assume now that $\displaystyle  I= \bigcup_{e\in \mathcal{U}(I)}W^1(e) = \bigcup_{e\in \mathcal{U}(I)} W_2 (e).$ To prove that $I$ is an inner ideal in $W$ we need to show that for any $a\in I$ and $b\in W,$ we have $\{a,b,a\}\in I$. Let $r_{_I} (a)\in I$ be the range tripotent of $a$ in the JBW$^*$-triple $I$. Clearly, $a\in I_2(r_{_I}(a)) \subseteq W_2(r_{_I}(a))$, by construction. It then follows from  Peirce arithmetic that $\{a,b,a\}\in  W_2(r_{_I}(a)) \subseteq I$ by our assumptions, which proves that $I$ is an inner ideal in $W$.
\end{proof}

The characterization of those norm-closed subtriples of a general JB$^*$-triple which are inner ideals is established next. 

\begin{proposition}\label{algebraic charac of inner ideals without dual}
Let $E$ be a real JB$^*$-triple, and let $I \subseteq E$ be a norm-closed subtriple. Then, the following conditions are equivalent:
\begin{enumerate}[$(i)$]
\item $I$ is an inner ideal in $E$.
\item $\displaystyle I =\bigcup_{a\in I} \left((E^{**})^1(r(a))\cap E\right)  = \bigcup_{a\in I} \left(E^{**}_2(r(a))\cap E\right)$.
\end{enumerate}
\end{proposition}

\begin{proof} Let $\overline{I}^{w^*}$ denote the weak$^*$-closure of $I$ in $E^{**}$. It follows from the separate weak$^*$-continuity of the triple product of $E^{**},$ and the natural embedding of $I$ inside its bidual space which identifies the latter with $\overline{I}^{w^*},$ that $I$ is a closed inner ideal of $E$ if, and only if, $\overline{I}^{w^*}$ is a weak$^*$-closed inner ideal of $E^{**}$. By  \Cref{algebraic charac of inner ideals} the latter statement is equivalent to  $$ \overline{I}^{w^*}    =\bigcup_{e\in \mathcal{U}\left(\overline{I}^{w^*}\right)} (E^{**})^1(e)  = \bigcup_{e\in \mathcal{U}\left(\overline{I}^{w^*}\right)} E^{**}_2(e) .$$ For each $a\in I$ we have $(E^{**})^1(r(a))\cap E  \subseteq   E^{**}_2(r(a)) \cap E\subseteq \overline{I}^{w^*}\cap I = I.$ We have therefore shown that $$\displaystyle \bigcup_{a\in I} \left((E^{**})^1(r(a))\cap E\right)  \subseteq \bigcup_{a\in I} \left(E^{**}_2(r(a))\cap E\right)\subseteq I.$$ Conversely, we simply observe that for each $a\in I,$ $a\in (E^{**})^1(r(a))\cap E  \subseteq   E^{**}_2(r(a)) \cap E.$
\end{proof}

Let $A$ be a real JB$^*$-algebra. As in the case of real C$^*$-algebras \cite[\S 5.2]{BingrenBook2003}, a functional $\varphi$ on ${A}$ is called \emph{positive} if it maps positive elements in $A$ to non-negative real numbers and $\varphi|_{A_{skew}} \equiv 0$. The symbol $(A^*)^+$ will stand for the set of all positive functional on $A$. It follows from this definition that  $(A^*)^+$ coincides with the set of all positive functionals on the JB-algebra $A_{sa}$. If $A$ is a real JBW$^*$-algebra with pedual $A_*$, we shall write $(A_*)^+$ or $A_*^+$ for the set of all positive normal functionals on $A$.\smallskip

It is not hard to see, via the characterization of inner ideals in complex JB$^*$-triples established by Edwards and Rüttimann in \cite{EdwRutti1992} and complexification, that every weak$^*$-closed inner ideal of a real JBW$^*$-triple is weak$^*$-Hahn--Banach smooth. The next theorem provides an algebraic characterization of those weak$^*$-closed subtriples of a real JBW$^*$-triple which are weak$^*$-Hahn--Banach smooth. As a consequence, we shall rediscover that every weak$^*$-closed inner ideal in a real JBW$^*$-triple is weak$^*$-Hahn--Banach smooth. However, contrary to what is known in the complex setting, the reciprocal implication does not necessarily hold.

\begin{theorem}\label{t characterization of wstarHahnBanach smooth subtriples} Let $I$ be a weak$^*$-closed subtriple of a real JBW$^*$-triple $W$. Then the following statements are equivalent:\begin{enumerate}[$(a)$]\item $I$ is weak$^*$-Hahn--Banach smooth, that is, every weak$^*$-continuous linear functional on $I$ admits a unique weak$^*$-continuous norm-preserving linear extension to $W$.
\item $\displaystyle I =\bigcup_{e\in \mathcal{U}(I)} W^1(e).$
\end{enumerate}
Consequently every weak$^*$-closed inner ideal of $W$ is weak$^*$-Hahn--Banach smooth.
\end{theorem}

\begin{proof} $(b) \Rightarrow(a)$ Suppose that $\displaystyle I= \bigcup_{e\in \mathcal{U}(I)} W^1(e)$ is a weak$^*$-closed subtriple in $W$. Take an arbitrary non-zero weak$^*$-continuous functional $\varphi \in I_{*}$. It is known that the set $\{e\in \mathcal{U}(I) : \|\varphi\| = \varphi (e) \}$ of all tripotents in $I$ supported at the functional $\varphi$ is non-empty (for example, the support tripotent of $\varphi$ in $I$ lies in this set, cf. \cite[Lemma 2.2]{Peralta2001little} and subsequent comments). Take a non-zero tripotent $e\in I$ satisfying $\varphi (e)= \|\varphi\|$. Clearly $e$ is a tripotent in $W$, and since $W^{1}(e) \subseteq I$ by hypothesis, we have $I^{1}(e) = W^{1}(e) \cap I = W^{1}(e)$. The projection $P^{1} (e): W\to W^{1}(e)$ is weak$^*$-continuous (see page~\pageref{ref Peirce projections and upper projections are weakstar conts}). Therefore, the functional $\Tilde{\varphi}:= \varphi \circ P^{1}(e) : W\to \mathbb{R}$ is weak$^*$-continuous with $\|\varphi\| = \varphi(e) =\Tilde{\varphi} (e)\leq \|\Tilde{\varphi}\| \leq\|\varphi\|$. \smallskip
	
We have therefore proved that every weak$^*$-continuous functional on $I$ admits a norm-preserving extension to a weak$^*$-continuous function on $W$. We have obtained this extension by making use of the hypothesis $(b)$. Actually a more general statement holds. Namely, by \cite[COROLLARY]{Bunce2001}, every weak$^*$-continuous functional on a weak$^*$-closed JB$^*$-subtriple of $W$ admits a norm-preserving extension to a weak$^*$-continuous functional on $W$. This result generalizes a celebrated result by S. Sakai in the setting of von Neumann algebras (see \cite[Proposition 1.24.5]{SakaiBook71}). \smallskip   
	
Concerning the uniqueness of the extension, suppose that $\psi\in W_*$ is another norm-preserving extension of $\varphi$. Let $e\in  I$ be the tripotent considered in the first paragraph. Since $ \|\psi\| = \|\varphi\|  = \varphi (e) = \psi(e),$ we deduce from  \cite[Lemma 2.7]{PeStacho2001} that $\psi = \psi\circ  P^{1} (e)= \psi|_{W^{1}(e)}\circ  P^{1} (e)$ on $W$. As before, $e\in I$ implies that $W^{1} (e) = I^{1} (e) $, and thus $\psi|_{W^{1}(e)} =\psi|_{I^{1}(e)} = \varphi|_{I^{1}(e)}.$ All together shows that $$\psi =  \psi|_{W^{1}(e)}\circ  P^{1} (e) = 
	\psi|_{I^{1}(e)}\circ  P^{1} (e)
	= \varphi|_{I^{1}(e)}\circ  P^{1} (e) = \Tilde{\varphi},$$ 
which concludes the proof of the uniqueness of the extension.\smallskip

$(a) \Rightarrow (b)$ Suppose now that every weak$^*$-continuous linear functional on $I$ has a unique weak$^*$-continuous norm-preserving extension to $W$. Our goal is to prove that $I$ satisfies the statement in $(b)$. We begin by establishing the following:\smallskip 

\noindent\emph{Claim 1: For every $e \in \mathcal{U}(I)$, every weak$^*$-continuous positive functional $\phi$  on the JBW-algebra $I^{1}(e)$ admits a unique weak$^*$-continuous norm-preserving Hahn--Banach extension $\Tilde{\phi}$ in $W^{1}(e)_{*}^{+}$.}\smallskip

In order to prove the claim fix $e\in \mathcal{U}(I)$, and $\phi \in (I^{1}(e))_{*}^{+}$. It is well known that $\|  \phi\| = \phi(e)$ (cf. \cite[Lemma 1.2.2]{HanchOlsenStormerBook}).\smallskip 

By observing that $e$ is a tripotent in the real JBW$^*$-triples $I$ and $W$, we can consider the projections $P_{_W}^{1}(e) : W \rightarrow W^{1}(e)$ and $P_{_I}^{1}(e): I \rightarrow I^{1}(e)$, which are both weak$^*$-continuous (cf. page~\pageref{ref Peirce projections and upper projections are weakstar conts}). We also know that $P^{1}_{_W} (e) \arrowvert_{I}  = P_{_I}^{1}(e)$ (see the formulae in page~\pageref{ref Peirce projections and upper projections are weakstar conts}). The functional $\varphi = \phi \circ P_{_I}^{1} (e) : I\to \mathbb{R}$ is weak$^*$-continuous on $I$ with $\|\varphi\| \leq \|\phi\| =\phi (e) = \varphi(e) \leq \|\varphi\|$. By the assumption on $I$,  there exists a (unique) weak$^*$-continuous functional $\Tilde{\varphi} \in W_*$ satisfying $\Tilde{\varphi}|_{_I} = {\varphi}$ and  $\| \Tilde{\varphi}\| = \| {\varphi} \| = \| \phi \|$. Let $\Tilde{\phi} = \Tilde{\varphi}|_{_{W^{1}(e)}} \in \left(W^{1} (e)\right)_*.$ By construction we have $$\|\phi \| = \| \Tilde{\varphi}\| \geq  \|\Tilde{\phi} \| \geq \Tilde{\phi} (e) = \phi (e) = \|\phi\|.$$ Since $e$ is the unit element of the JBW-algebra $W^1 (e)$ we deduce that $\Tilde{\phi}$ is a positive normal functional on $W^{1} (e)$ (cf. \cite[Lemma 1.2.2]{HanchOlsenStormerBook}). Clearly $\Tilde{\phi}$ is a positive norm-preserving extension of $\phi$ to $W^{1} (e)$. \smallskip

For the uniqueness in the claim, consider another functional $\psi \in (W^{1}(e))_{*}^{+}$ such that $\psi\arrowvert_{I^{1}(e)} = \phi$ and $\|\phi \| = \|\psi\| =\psi (e)$. As before, by hypotheses, the functional $\rho=\psi \circ P^{1}_{_I} (e) : I\to \mathbb{R}$ admits a unique norm-preserving extension to a functional $\Tilde{\rho}\in W_*$. Since $$\Tilde{\rho}|_{_I} = \psi \circ P^{1}_{_I} (e) = \psi\arrowvert_{I^{1}(e)} \circ P^{1}_{_I} (e) = \phi \circ P^{1}_{_I} (e) = \varphi,$$ with $\|\Tilde{\rho} \| = \| \rho\| = \|\psi\| = \|\phi\|,$ it follows from the uniqueness of the norm-preserving weak$^*$-continuous extension of $\varphi$ to $W$ that $\Tilde{\rho} = \Tilde{\varphi}$. Consequently, $\psi = \Tilde{\rho}|_{_{W^{1}(e)}} = \Tilde{\varphi} |_{_{W^{1}(e)}} = \Tilde{\phi}$, which concludes the proof of the claim.\smallskip

We observe now that $I^{1} (e)$ is a weak$^*$-closed subtriple of $W^{1} (e)$, and hence a weak$^*$-closed Jordan subalgebra of the latter. We also know that $I^{1} (e)$ enjoys the property stated in Claim 1. Therefore, \cite[Lemma 2.1]{EdwRutti1992} guarantees that  $I^{1}(e)$ is an inner ideal in $W^{1}(e)$. In particular, $W^{1}(e) = I^{1}(e),$ and thus 
$$ \bigcup_{e\in \mathcal{U}(I)} W^1(e) = \bigcup_{e\in \mathcal{U}(I)} I^1(e) \subseteq I.$$ We finally observe that for each $a\in I$ we have $a =\{r_{_I}(a), a, r_{_I}(a)\}\in I^1 (r_{_I}(a)) = W^{1} (r_{_I}(a))$, which concludes the proof of $(a)\Rightarrow (b)$.\smallskip

The final statement is a consequence of \Cref{algebraic charac of inner ideals} and the previous equivalence. 
\end{proof}

\begin{remark}\label{remark minimal tripotents in Hahn--Banach smooth subtriples}
Suppose $I$ is a weak$^*$-closed JB$^*$-subtriple of a real JBW$^*$-triple $W$. Minimal tripotents in $I$ need not be minimal in $W$. If $I$ is an inner ideal in $W$, it is not hard to check that $\mathcal{U}_{min}(I) \subseteq \mathcal{U}_{min}(E)$, since for each $e\in \mathcal{U}_{min}(I)$, we have $E^{1} (e) = I^{1} (e) = \mathbb{R} e$. If we only assume that $I$ is weak$^*$--Hahn--Banach smooth, the containing  $\mathcal{U}_{min}(I) \subseteq \mathcal{U}_{min}(E)$ is also true thanks to Theorem~\ref{t characterization of wstarHahnBanach smooth subtriples}.
\end{remark}

We state next a couple of corollaries characterising norm closed (resp., weak$^*$-closed) inner ideals in a real JB$^*$-triple (resp., a real JBW$^*$-triple) among subtriples which are Hahn Banach smooth (resp., weak$^*$-Hahn Banach smooth).

\begin{corollary}\label{c characterization of weak* inner ideals} Let $I$ be a weak$^*$-closed subtriple of a real JBW$^*$-triple $W$. Then the following statements are equivalent:\begin{enumerate}[$(a)$]\item $I$ is an inner ideal.
\item $I$ is weak$^*$-Hahn--Banach smooth and for each tripotent $e\in I$ we have $W^{-1}(e)\subseteq I.$
\end{enumerate}
\end{corollary}

\begin{proof} $(a)\Rightarrow (b)$ is clear from Theorem~\ref{t characterization of wstarHahnBanach smooth subtriples}. For the reciprocal implication observe that if $I$ satisfies $(b)$, it follows from \Cref{t characterization of wstarHahnBanach smooth subtriples} that $\displaystyle I =\bigcup_{e\in \mathcal{U}(I)} W^1(e).$ However, the extra assumptions on $I$ assure that $$  \bigcup_{e\in \mathcal{U}(I)} W_2(e) =\bigcup_{e\in \mathcal{U}(I)} W^{1}(e) \oplus W^{-1} (e) \subseteq I \subseteq  \bigcup_{e\in \mathcal{U}(I)} W_2(e) .$$ \Cref{algebraic charac of inner ideals} implies that $I$ is an inner ideal. 
\end{proof}

\begin{corollary}\label{c characterization of closed inner ideals} Let $I$ be a closed subtriple of a real JB$^*$-triple $E$. Then the following statements are equivalent:\begin{enumerate}[$(a)$]\item $I$ is an inner ideal.
\item $I$ is Hahn--Banach smooth and for each norm-one element $a\in I$ we have $(E^{**})^{-1}(r(a))\cap E\subseteq I.$
\end{enumerate}
\end{corollary}

\begin{proof} $(a) \Rightarrow (b)$ It can be easily deduced from the separate weak$^*$-continuity of triple product (cf. \cite{MartinezPe2000separate}) and the weak$^*$-density of $I$ inside $I^{**}$ that $I$ is a closed inner ideal in $E$ if, and only if, $I^{**}$ is a weak$^*$-closed inner ideal in $E^{**}$, and by \Cref{c characterization of weak* inner ideals}, the latter holds if, and only if, $I^{**}$ is weak$^*$-Hahn--Banach smooth in $E^{**}$ and for each tripotent $e\in I^{**}$ we have $(E^{**})^{-1}(e)\subseteq I^{**}.$ It clearly follows that $I$ is Hahn--Banach smooth and for each norm-one element $a\in I$ we have $(E^{**})^{-1}(r(a))\cap E\subseteq I^{**}\cap E =  I$ (cf. the Bipolar theorem).\smallskip

\noindent $(b) \Rightarrow (a)$ As above $I$ is Hahn--Banach smooth in $E$ if, and only if, $I^{**}$ is weak$^*$-Hahn--Banach smooth in $E^{**}$. \Cref{t characterization of wstarHahnBanach smooth subtriples} implies that $\displaystyle I^{**} = \cup_{e\in \mathcal{U}(I^{**
})} \left(E^{**}\right)^{1} (e).$ It then follows that $\left(E^{**}\right)^{1} (r(a))\cap E \subseteq I,$ for all norm-one element $a\in I.$ The remaining assumption on $I$ assures that $$\left(E^{**}\right)^{-1} (r(a))\cap E \subseteq I,$$ and thus $\left(E^{**}\right)_{2} (r(a))\cap E \subseteq I$ for all norm-one element $a\in I.$  \smallskip

Finally, for each norm-one element $a\in I$ and every $x\in E,$ by Peirce arithmetic we have $\{a,x,a\}\in (E^{**})_2 (r(a)) \cap E \subseteq I.$ It is then clear that $I$ is an inner ideal in $E$.
\end{proof}

It should be point out that in the setting of complex JB$^*$-triples, the second condition in \Cref{c characterization of weak* inner ideals}$(b)$ can be relaxed. Indeed, if $\mathcal{I}$ is a weak$^*$-closed subtriple of a JBW$^*$-triple $\mathcal{W}$, for each tripotent $e\in \mathcal{I}$ we have $\mathcal{W}^{-1} (e) = i \mathcal{W}^{1} (e) \subseteq \mathcal{I}.$ However, the next example shows that weak$^*$-Hahn--Banach smoothness is not, in general, equivalent to being an inner ideal in the real setting.

\begin{example}\label{example w*HahnBanach smoothness is not enough} Let $W = \mathbb{C}$ regarded as a real JBW$^*$-triple, and let $I = \mathbb{R}$. Clearly, $\mathcal{U} (I) = \{0,1,-1\}$, $\displaystyle I =\bigcup_{e\in \mathcal{U}(I)} W^1(e)$, and hence $I$ is weak$^*$-Hahn--Banach smooth, but $I$ is not an inner ideal of $W$. There is a simple explanation for this counterexample. For a complex Banach space $X$, whose underlying real Banach space is denoted by $X_{_{\mathbb{R}}}$, the dual space of $X_{_\mathbb{R}}$ is not big enough. We can provide more natural examples. Take a complex JBW$^*$-triple $\mathcal{W}$ and a minimal tripotent $e\in \mathcal{W}$. By setting $I =\mathcal{W}^{1} (e) = \mathbb{R} e,$ we get a weak$^*$-closed subtriple which is not an inner ideal since $\mathcal{W}^{-1} (e) \cap I = \{0\}.$  It is clear that $I$ is weak$^*$-Hahn--Banach smooth in $\mathcal{W}_{_\mathbb{R}}$.\smallskip
	
Let $E = c_0$ denote the JB$^*$-triple of all null sequences fo complex numbers, and let $I$ denote the real JB$^*$-subtriple of $E$ of all real sequences whose components in a determined set $\emptyset \subsetneq N_0\subsetneq \mathbb{N}$ vanish. It is not hard to check, via the arguments in the previous paragraph, that $I$ is Hahn-Banach smooth but it is not an inner ideal of $E$. 	
\end{example}

After the previous example, the challenge now consists in determining whether, under some extra hypotheses, in some real JBW$^*$-subtriples, weak$^*$-Hahn--Banach smoothness is enough to guarantee that a weak$^*$-closed subtriple is an inner ideal. We shall get some positive and surprising answers in Sections~\ref{sec: reduced Cartan factors}, \ref{sec: complex Cartan factors} and \ref{sec: conclusions} below.\smallskip

The characterization of closed inner ideals in a real JB$^*$-triple $E$ given in \Cref{c characterization of closed inner ideals} benefits from the characterization in \Cref{c characterization of weak* inner ideals} because the weak$^*$-closure of each inner ideal in $E$ is an inner ideal in $E^{**}$. However, the extension of \Cref{t characterization of wstarHahnBanach smooth subtriples} to the case of closed subtriples which are Hahn--Banach smooth will require a more elaborated machinery. Note that if $I$ is a norm closed real JB$^*$-subtriple of a real JB$^*$-triple, $I$ being Hahn--Banach smooth in $E$ is equivalent to $\overline{I}^{w^*}$ being weak$^*$ Hahn-Banach smooth in $E^{**}$. However, a direct application of Theorem~\ref{t characterization of wstarHahnBanach smooth subtriples} requires us to deal with all tripotents in $\overline{I}^{w^*}$ not only with those associated with range tripotents of elements in $I$.\smallskip

We begin with a quantitative version of \cite[Lemma 2.7]{PeStacho2001}. Let us first recall some definitions. Suppose  $\mathcal{U}$ is an ultrafilter on a given index set $\Lambda$, and let $(X_i)_{i\in \Lambda}$ be a family of Banach spaces. The set $$N_{\mathcal{U}} := \left\{ (x_i)_i\in \bigoplus_{i\in \Lambda}
^{\ell_\infty} X_i : \lim_{\mathcal{U}} \|x_i\| =0 \right\}$$ is a closed subspace of the Banach space $\bigoplus_{i\in \Lambda}
^{\ell_\infty} X_i,$ and the quotient space $(\bigoplus_{i\in \Lambda}
^{\ell_\infty} X_i) / N_{\mathcal{U}}$ is known as the $\mathcal{U}$-ultraproduct of the $X_i$'s and is denoted by $(X_i)_{\mathcal{U}}$. Elements in this quotient space will be denoted by $[x_i]_{\mathcal{U}}$.  In case that $X_i=X$ for all $i$, we write $(X)_{\mathcal{U}}$ for the corresponding ultraproduct and we call it the ultrapower of $X$. If every $X_i$ is a JB$^*$-triple, the space $\bigoplus_{i\in \Lambda}
^{\ell_\infty} X_i$ is a JB$^*$-triple with respect to the natural point-wise triple product \cite[page 523]{Kaup1983riemann}, and the subspace $ N_{\mathcal{U}}$ is a closed ideal of it. This shows that the ultraproduct $(X_i)_{\mathcal{U}}$ is a JB$^*$-triple (cf. \cite[Corollary 10]{Dineen1986} or \cite[page 157]{PeraltaPfitzner2016}), and consequently the ultraproduct of a family of real JB$^*$-triples is a real JB$^*$-triple.

\begin{lemma}\label{l: quatitative version of PeSta P1} For each $\varepsilon >0$ there exists a positive number $\delta$ satisfying the following property: for each real JB$^*$-triple $E$, every tripotent $e\in E$ and every $\phi\in E^{*}$ satisfying $\phi (e) > \|\phi\| -\delta$, we have $\|\phi - \phi \circ P^{1} (e)\| <\varepsilon.$
\end{lemma}

\begin{proof} We can clearly reduce to the case in which $\|\phi\| =1$. Arguing by contradiction, we assume the existence of a positive $\varepsilon$ such that for each positive $\delta$ there exists a real JB$^*$-triple $E_{\delta}$, a tripotent $e_{\delta}\in E_{\delta}$ and a norm-one functional $\phi_{\delta}\in E_{\delta}^*$ satisfying $\phi_{\delta} (e_{\delta}) > 1-\delta$ but  $\|\phi_{\delta} - \phi_{\delta} \circ P^{1} (e_{\delta})\| \geq \varepsilon.$ Let $\mathcal{U}$ be a free ultrafilter over $\mathbb{N}$, and consider the ultraproduct $(E_n)_{\mathcal{U}}$, where $E_n$ is the real JB$^*$-triple whose existence has been assumed above for $\delta = \frac{1}{n}$ ($n\in \mathbb{N}$). We have already commented that $(E_n)_{\mathcal{U}}$ is a real JB$^*$-triple. It is easy to check that the element $[e_n]_{\mathcal{U}}$ is a tripotent in $(E_n)_{\mathcal{U}}$.\smallskip

It is a well known result that the element $[\phi_{n}]_{\mathcal{U}}\in (E_n^*)_{\mathcal{U}}$ defines a bounded linear functional on $(E_n)_{\mathcal{U}}$ given by $[\phi_{n}]_{\mathcal{U}} ([x_{n}]_{\mathcal{U}}) := \lim_{\mathcal{U}} \phi_n (x_n)$ whose norm is precisely $\| [\phi_{n}]_{\mathcal{U}} \|$, so $\| [\phi_{n}]_{\mathcal{U}} \|=1$ in this case (see \cite[\S 7]{Heinrich1980}). It is easy to check from the definition of the triple product on $(E_n)_{\mathcal{U}}$ that $$P^{1} ([e_n]_{\mathcal{U}}) ([x_{n}]_{\mathcal{U}}) = [P^{1} (e_n) (x_{n})]_{\mathcal{U}}, \hbox{for all } [x_{n}]_{\mathcal{U}} \in (E_n)_{\mathcal{U}}.$$ Finally, it follows from the assumptions that $$[\phi_{n}]_{\mathcal{U}} ([e_{n}]_{\mathcal{U}}) := \lim_{\mathcal{U}} \phi_n (e_n) = 1,$$ and $$\begin{aligned}
 \left\| [\phi_{n}]_{\mathcal{U}} - [\phi_{n}]_{\mathcal{U}} \circ P^{1} ([e_n]_{\mathcal{U}}) \right\| &= \left\| [\phi_{n} - \phi_{n} \circ P^{1} (e_n)]_{\mathcal{U}}) \right\| \\
 &= \lim_{\mathcal{U}} \| \phi_{n} - \phi_{n} \circ P^{1} (e_n)\| \geq \varepsilon>0,   
\end{aligned}$$ which contradicts \cite[Lemma 2.7]{PeStacho2001}.
\end{proof}

The desired characterisation of norm closed subtriples which are Hahn--Banach smooth reads as follows:

\begin{theorem}\label{t characterization of closed subtriples with HB smoothness} Let $I$ be a closed subtriple of a real JB$^*$-triple $E$. Then the following statements are equivalent:\begin{enumerate}[$(a)$]\item $I$ is Hahn--Banach smooth.
\item $\displaystyle I =\bigcup_{a\in I, \|a\| =1} (E^{**})^{1}(r(a))\cap E.$
\end{enumerate}
\end{theorem}

\begin{proof} $(a)\Rightarrow (b)$ As commented in previous results (cf. \Cref{c characterization of closed inner ideals}), $I^{**}$ is a weak$^*$-closed subtriple of $E^{**}$ which is weak$^*$-Hahn--Banach smooth. Thus, by \Cref{t characterization of wstarHahnBanach smooth subtriples}, $\displaystyle I^{**} = \bigcup_{e\in \mathcal{U}(I^{**
})} (E^{**})^{1}(e).$ We deduce that $I =  E\cap I^{**} = \bigcup_{e\in \mathcal{U}(I^{**})} \left((E^{**})^{1}(e)\cap E\right)$, in particular, $(E^{**})^{1}(r(a))\cap E \subseteq  I$, for all norm-one element $a\in I$. Reciprocally, each norm-one element $a\in I$ lies in $(E^{**})^{1}(r(a))\cap E \subseteq  I$.\smallskip

$(b)\Rightarrow (a)$ This implication gives rise to some additional difficulties and is more than a mere application of \Cref{t characterization of wstarHahnBanach smooth subtriples}. Suppose that $(b)$ holds. It is enough to show the uniqueness in the Hahn--Banach smooth property.
Let $\phi$ be a norm-one functional in $I^*$, and let $\tilde{\phi},\psi\in E^*$ be two norm-preserving extensions of $\phi.$\smallskip 

Fix an arbitrary $\varepsilon>0,$ and the corresponding $\delta>0$ given by \Cref{l: quatitative version of PeSta P1} for $\frac{\varepsilon}{2}$. Since $\|\phi\|=1$, we can find a norm-one element $a_{_\delta}\in I$ such that $\phi (a_{_\delta}) > 1-\frac{\delta}{2}.$ The functionals $\tilde{\phi}|_{_{(E^{**})^{1}(r(a_{_\delta}))}},$ $\psi|_{_{(E^{**})^{1}(r(a_{_\delta}))}}$ can be written in the form $$\tilde{\phi}|_{ _{(E^{**})^{1}(r(a_{_\delta}))}} = \tilde{\phi}^+ -\tilde{\phi}^{-},\ \ \ \psi|_{_{(E^{**})^{1}(r(a_{_\delta}))}} = \psi^+ -\psi^-,$$ where $\tilde{\phi}^+,\tilde{\phi}^{-}, \psi^+,$ and $\psi^-$ are positive normal functionals on the real JBW$^*$-algebra $(E^{**})^{1}(r(a_{_\delta}))$ with $$\begin{aligned}
1\geq& \|\tilde{\phi}|_{_{(E^{**})^{1}(r(a_{_\delta}))}}\| = \|\tilde{\phi}^+\| + \|\tilde{\phi}^{-}\|, \hbox{ and }\\ 1\geq& \|\psi|_{_{(E^{**})^{1}(r(a_{_\delta}))}}\| = \|\psi^+\| +\|\psi^-\| 
\end{aligned}$$ (cf. \cite[Proposition 4.5.3]{HanchOlsenStormerBook}). Since $$ 1-\frac{\delta}{2} < \phi (a_{_\delta}) = \tilde{\phi}|_{ _{(E^{**})^{1}(r(a_{_\delta}))}} = \tilde{\phi}^+ (a_{_\delta}) -\tilde{\phi}^{-} (a_{_\delta}),$$ we deduce that $$\| \tilde{\phi}^+ \| \geq \tilde{\phi}^+ (a_{_\delta}) > 1-\frac{\delta}{2}, \hbox{ and } \|\tilde{\phi}^{-}\| <\frac{\delta}{2},$$ and thus $\tilde{\phi} (r(a_{_\delta})) = \tilde{\phi}|_{_{(E^{**})^{1}(r(a_{_\delta}))}} (r(a_{_\delta})) = \tilde{\phi}^+ (r(a_{_\delta})) -\tilde{\phi}^{-} (r(a_{_\delta})) > 1-\delta.$ \Cref{l: quatitative version of PeSta P1} implies that $\| \tilde{\phi} - \tilde{\phi} \circ P^{1} (r(a_{_\delta})) \| <\frac{\varepsilon}{2}.$ We can similarly obtain that $\psi (r(a_{_\delta})) > 1-\frac{\delta}{2}$ and $\| {\psi} - {\psi} \circ P^{1} (r(a_{_\delta})) \| <\frac{\varepsilon}{2}.$\smallskip

Now, having in mind that $a_{_\delta}$ is a norm-one elemen in $I$, by assumptions $(E^{**})^{1}(r(a_{_\delta}))\cap E \subseteq I$, the JB-algebra $E(a_{_\delta})$ is contained in $(E^{**})^{1}(r(a_{_\delta}))\cap E$ hence in $I$ and is weak$^*$-dense in $(E^{**})^{1}(r(a_{_\delta}))$ (see page~\pageref{eq density properties of the self-adjoint part}), and the functionals $\tilde{\phi}$ and $\psi$ are weak$^*$-continuous on $E^{**}$ with $\tilde{\phi}|_{I} = \phi = \psi|_{I}$, we obtain $\tilde{\phi}|_{(E^{**})^{1}(r(a_{_\delta}))} =  \psi|_{(E^{**})^{1}(r(a_{_\delta}))}$. \smallskip

Finally, summarizing the above conclusions, we arrive to  
$$\begin{aligned} \| \tilde{\phi}- \psi\| &\leq \| \tilde{\phi} - \tilde{\phi}\circ P^{1} (r(a_{_\delta}))\| + \| {\psi} - {\psi}\circ P^{1} (r(a_{_\delta}))\|   \\ &+  \| \tilde{\phi}|_{(E^{**})^{1}(r(a_{_\delta}))} \circ P^{1} (r(a_{_\delta})) - {\psi}|_{(E^{**})^{1}(r(a_{_\delta}))}\circ P^{1} (r(a_{_\delta}))\| <\varepsilon.
\end{aligned}$$ The arbitrariness of $\varepsilon>0$ implies that $\tilde{\phi}= \psi,$ and hence $I$ is Hahn--Banach smooth.
\end{proof}

We have seen in Example~\ref{example w*HahnBanach smoothness is not enough} that there exists JB$^*$-subtriples which are Hahn-Banach smooth without being inner ideals. However there is an algebraic notion which is valid to characterise Hahn-Banach smoothness of JB$^*$-subtriples. Let $F$ be a real JB$^*$-subtriple of a real JB$^*$-triple $E$. Following the ideas in \cite[\S 4]{FerPe2010inner} and inspired by the usual notion of hereditary C$^*$-subalgebra, we say that $F$ is a \emph{hereditary subtriple} of $E$ if for each $a\in F$, the real JB$^*$-algebra $F(a)$ behaves hereditarily as a real JB$^*$-subalgebra of $E(a)$, that is, if $b\in F(a)$ and $c\in E(a)$ are two positive elements with $c\leq b$ in $E(a)$, we have $c\in F$. Note that the inner ideal $F(a)$ is a JB$^*$-subalgebra of $E(a)$ (cf. page~\pageref{eq density properties of the self-adjoint part}).

\begin{theorem}\label{t characterization of closed subtriples with HB smoothness as hereditary subtriples} Let $I$ be a closed subtriple of a real JB$^*$-triple $E$. Then the following statements are equivalent:\begin{enumerate}[$(a)$]\item $I$ is Hahn--Banach smooth.
\item $I$ is a hereditary subtriple of $E$.
	\end{enumerate}
\end{theorem}

\begin{proof} $(a)\Rightarrow (b)$ Suppose $I$ is Hahn-Banach smooth. Theorem~\ref{t characterization of closed subtriples with HB smoothness} assures that $\displaystyle I =\bigcup_{a\in I, \|a\| =1} (E^{**})^{1}(r(a))\cap E.$  Let us fix a norm-one element $a\in I$ and positive elements $b\in I(a)$, $c\in E(a)$ with $c\leq b$. Observe that $c\in (E^{**})^{1}(r(a))\cap E \subseteq I$.\smallskip

$(b)\Rightarrow (a)$ Take now a norm-one element $a\in I$, and $x\in (E^{**})^{1}(r(a))\cap E$. Working in the JBW-algebra $(E^{**})^{1}(r(a))$ we can write $x = x^{+} - x^{-}$, where $x^+,x^-$ are two orthogonal positive elements in $(E^{**})^{1}(r(a))\cap E.$ It is known that the sequence $(a^{[1/(2n-1)]})_{n}$ converges to $r(a)\in I^{**} = \overline{I}^{w^*}$ in the strong$^*$-topology of $I^{**}$ (and of $E^{**}$) as defined in \cite[\S 4]{PeRod2001Groth} (see also \cite[Corollary]{Bunce2001}). Note that $0 \leq U_{_{a^{[1/(2n-1)]}}} (x^{+}), U_{_{a^{[1/(2n-1)]}}} (x^{-})\leq U_{_{a^{[1/(2n-1)]}}} (r(a)),$ for all natural $n$, where $U_{_{a^{[1/(2n-1)]}}} (r(a))\in I$. The assumptions on $I$ imply that $U_{_{a^{[1/(2n-1)]}}} (x^{+})$ and $U_{_{a^{[1/(2n-1)]}}} (x^{-})$ belong to $I(a)$ for all natural $n$. Having in mind that the triple product of $E^{**}$ is jointly strong$^*$ continuous on bounded sets, together with the fact that the strong$^*$-topology is stronger that the weak$^*$-topology (cf. \cite[Theorem 9 and Corollary 9]{PeRod2001Groth}), by taking strong$^*$-limits in $n$ in the above expressions, we get $x^{+} = U_{_{r(a)}} (x^{+})$ and $x^{-} = U_{_{r(a)}} (x^{-})$ belong to $\overline{I(a)}^{w^*}\subseteq  \overline{I}^{w^*} = I^{**}.$ Since $x^{+},x^{-}\in E(a)\cap \overline{I(a)}^{w^*} = I(a).$ Consequently, $x = x^+- x^{-}\in I(a) \subseteq I.$ The arbitrariness of $x \in (E^{**})^{1}(r(a))\cap E$ proves that $(E^{**})^{1}(r(a))\cap E \subseteq I$, for all norm-one element $a\in I$, which in view of Theorem~\ref{t characterization of closed subtriples with HB smoothness} guarantees that $I$ is Hahn-Banach smooth.
\end{proof}

When in the proof of Theorem~\ref{t characterization of closed subtriples with HB smoothness as hereditary subtriples}, Theorem~\ref{t characterization of wstarHahnBanach smooth subtriples} replaces Theorem~\ref{t characterization of closed subtriples with HB smoothness} we get the following result. 

\begin{theorem}\label{t characterization of closed subtriples with weak*HB smoothness as hereditary subtriples} Let $I$ be a weak$^*$-closed subtriple of a real JBW$^*$-triple $W$. Then the following statements are equivalent:\begin{enumerate}[$(a)$]\item $I$ is weak$^*$-Hahn--Banach smooth.
\item $I$ is a hereditary subtriple of $W$.
\end{enumerate}
\end{theorem}

In the setting of complex JB$^*$-triples, hereditary JB$^*$-subtriples and inner ideals define the same objects. However in the case of real JB$^*$-triples, the collection of all inner ideals is strictly smaller (see Example~\ref{example w*HahnBanach smoothness is not enough}).\smallskip

Norm-preserving extensions of weak$^*$ continuous functionals on a weak$^*$-closed subspace $Y$ of a dual Banach space $X^*$ are not, in general, available (we can simply consider a norm-closed subspace $M$ of $X$ for which the quotient norm in $X/M$ is not attained at some point $\phi+M\in X/M$ and take $Y = M^{\circ} \cong \left(X/M\right)^*$, the polar of $M$ in $X^*$).\smallskip

Let $X$ be a (real or complex) Banach space and $M$ be a closed subspace of it. We say that $M$ \emph{has the Haar property in $X$} if for each $x \in X$ there exists a unique $y \in M$ such that $\| x-y\| = dist (x, M)$ where $ dist (x, M) := \inf \{ \| x-z\|: z\in M\}$. Denote $\mathcal{B}_{X}$ and $\mathcal{S}_{X}$ as the closed unit ball and unit sphere of $X,$ respectively. In the next proposition we revisit \cite[Theorem 1.1]{Phelps1960}.

\begin{proposition}
Let $X$ be a Banach space, and let $M$ be a weak$^*$-closed subspaces of $X^*$ satisfying that each weak$^*$-continuous functional on $M$ admits a norm-preserving extension to a weak$^*$-continuous functional on $X^*$. Then $M$ is weak$^*$-Hahn--Banach smooth (in $X^*$) if, and only if, $M_{\circ}$, the prepolar of $M$ in $X$, satisfies the Haar property in $X$.
\end{proposition}

\begin{proof} The proof in \cite[Theorem 1.1]{Phelps1960}, or even a simplified version of it, is also valid in this case. It is included here for completeness reasons. We claim that for each $\phi \in X$, $dist(\phi, M_{\circ}) = \| \phi\arrowvert_{M} \|$. Actually, given $\psi \in M_{\circ}$, we have $$ \| \phi\arrowvert _{M} \| = \sup \{\|(\phi-\psi)(m)\|: m\in M, \|m\|\leq 1\} \leq \|\phi-\psi\|,$$ which implies that $ \| \phi\arrowvert_{M} \| \leq dist(\phi, M_{\circ})$. Now, we apply the hypothesis on $M$ to find a functional $\varphi \in X$, which is a  norm-preserving weak$^*$-continuous extension of $\phi|_{M}$. Under these conditions $\phi-\varphi \in M_{\circ}$ and furthermore
	\begin{align*}
		\|\phi|_{M}\| &= \|\varphi\| = \sup\{ \|\varphi(x)\| : x\in X^*, \|x\|\leq 1\} \\
		&= \sup\{ \|(\phi-(\phi-\varphi)) (x)\| : x\in X^*, \|x\|\leq 1\} \\
		&= \| \phi-(\phi-\varphi)\| \geq dist(\phi,M_{\circ}),
	\end{align*} which finishes the proof of the claim.\smallskip 
	
If $M_{\circ}$ does not satisfy the Haar property in $X$, then there exist $\varphi\in X$ and $0\neq \phi\in M_{\circ}$ such that $$\|\varphi\arrowvert_{M} \| = dist(\varphi, M_{\circ}) = 1 = \|\varphi -\phi\|.$$ This implies that $\varphi$ and $\varphi-\phi$ are two different norm-preserving extensions of $\varphi\arrowvert_{M},$ which implies that $M$ is not weak$^*$-Hahn--Banach smooth.\smallskip

Reciprocally, if $M$ is not weak$^*$-Hahn--Banach smooth, there exist a norm-one functional $\phi \in {M_{*}}$ admitting two different norm-preserving weak$^*$-continuous extensions $\varphi,\psi \in X$. In such a case $0 \neq \varphi-\psi \in M_{\circ}$ and 
\begin{align*}
	1= \|\varphi\| &= \|\psi\| = \|\varphi-(\varphi-\psi)\| \\
	&\geq dist(\varphi,M_{\circ})= \|\varphi\arrowvert_{M} \| = \|\phi\|=1,
\end{align*} witnessing that $M_{\circ}$ does not satisfy the Haar property in $X$.
\end{proof}

By recalling that every weak$^*$-closed subtriple of a real JBW$^*$-triple $E$ satisfies the hypothesis in the previous proposition we get the following corollary.

\begin{corollary}\label{c to Phelps theorem for weakstar topology} Let $I$ be a weak$^*$-closed subtriple of a real JBW$^*$-triple $W$. Then $I$ is weak$^*$-Hahn--Banach smooth if, and only if, $I_{\circ}$, the prepolar of $I$, satisfies the Haar property in $W_{*}$.
\end{corollary}

\section{Inner ideals in reduced atomic real JBW$^*$-triples}\label{sec: reduced Cartan factors}

Recall that a real JB$^*$-triple $E$ is called \emph{reduced} if $E^{2} (e) = \mathbb{R}e$ (i.e. $E^{-1}(e) = \{0 \}$) for every minimal tripotent $e \in E$ (cf. \cite[11.9]{Loos1977bounded}). Every complex Cartan factor, regarded as a real JB$^*$-triple, is non-reduced. However, there is a wide list of real Cartan factors which are reduced real JB$^*$-triples. According to the terminology coined by Kaup in \cite{Kaup1997}, \emph{real Cartan factors} as defined as the real forms of (complex) Cartan factors and they can be classified into 12 different types (including 8 classical and 4 exceptional types). We introduce now all real Cartan factors. Let $\mathcal{H}$ be a complex Hilbert space of dimension $n$, $X,Y$ real Hilbert spaces of dimensions $n,m,$ respectively, and $P,Q$ Hilbert spaces of dimensions $p$ and $q$ over the quaternion field $\mathbb{H}$, respectively. The real Cartan factors are defined as follows:\vspace{1mm}

\begin{enumerate}[$(i)$]
\item $I_{n,m}^{\mathbb{R}} := \mathcal{L}(X,Y) \ (m\geq n\geq 2),$\vspace{1mm}
\item  $I_{2p,2q}^{\mathbb{H}}:=\mathcal{L}(P,Q) \ (q\geq p\geq 2),$\vspace{1mm}
\item  $I_{n,n}^{\mathbb{C}}:= \{z\in \mathcal{L}(\mathcal{H}): z^{*} = z\} \ (n\geq 2),$\vspace{1mm}
\item  $II_{n}^{\mathbb{R}}:= \{x \in \mathcal{L}(X): x^* = -x \}\ (n\geq 4)$,\vspace{1mm}
\item  $II_{2p}^{\mathbb{H}}:= \{ w \in \mathcal{L}(P): w^*=w \} \ (p\geq 2),$\vspace{1mm}
\item $III_{n}^{\mathbb{R}} := \{ x \in \mathcal{L}(X) : x^* = x\} \ (n\geq 2),$ \vspace{1mm}
\item $III_{2p}^{\mathbb{H}}:= \{ w \in \mathcal{L}(P): w^* = -w \} \ (p\geq 2),$\vspace{1mm}
\item (\emph{Real spin factor}) $IV_{n}^{r,s} := X_{1} \mathop{\oplus}\limits^{\ell_{1}} X_{2},$ where $X_{1}$ and $X_{2}$ are closed orthogonal linear subspaces of $X$ of dimension $r$ and $s,$ respectively (with $r+s\geq 3$, $r\geq s \geq 1$). The triple product of $E$ is defined by $$\{x,y,z \}  = \langle x|y\rangle z + \langle z| y\rangle x - \langle x| \bar{z}\rangle  \bar{y},$$
where $\langle \cdot| \cdot\rangle$ is the inner product of $X$ and the involution $\bar{\cdot}$ on $E$ is given by $\bar{x} := (x_{1}, -x_{2})$ for every $x = (x_{1}, x_{2}) \in E$. \vspace{1mm}
\end{enumerate}
$
\begin{array}{ll}
	{(ix) }\ V^{\mathbb{O}^{rs}} := M_{1,2}(\mathbb{O}^{rs}), & (x) \ V^{\mathbb{O}^{rd}}: = M_{1,2}(\mathbb{O}^{rd}), \\[8pt]
	{(xi) }\ VI^{\mathbb{O}^{rs}} := H_{3}(\mathbb{O}^{rs}), & {(xii) }\ VI^{\mathbb{O}^{rd}} := H_{3}(\mathbb{O}^{rd}),
\end{array}
$ \\[8pt] 
where $\mathbb{O}^{rs}$ is the real split Cayley algebra over $\mathbb{R}$ and $\mathbb{O}^{rd}$ is the real division Cayley algebra (also known as the algebra of real division octonions). The first 8 types are the classical types, while the real Cartan factors of types $(ix)$--$(xii)$ are called \emph{exceptional real Cartan factors}.\smallskip

According to \cite[11.9]{Loos1977bounded} and \cite[table 1]{Kaup1997}, the unique non-reduced real Cartan factors are $IV_{n}^{n,0}, V^{\mathbb{O}^{rd}}, I_{2p,2q}^{\mathbb{H}}$ and $III_{2p}^{\mathbb{H}}$ (corresponding to those factors with $z =1$ in the mentioned table). In this section, we shall prove that for each reduced atomic real JBW$^*$-triple $W$ (i.e. a real JBW$^*$-triple which is the direct sum of reduced real Cartan factors), every weak$^*$-closed subtriple which is weak$^*$-Hahn--Banach smooth is automatically an inner ideal.\smallskip

Let $E$ be a real JB$^*$-triple. A subset $S \subseteq E$ is called \emph{orthogonal} if $0 \notin S$ and $a \perp b$ for all $a,b \in F$. We say that \emph{$E$ has rank $m$} if $m$ is the minimal cardinal number such that $card(S) \leq m$ for every orthogonal subset $S$ in $E$ (cf. \cite[pp. 196]{Kaup1997}). A tripotent $u \in E$ has \emph{rank $m$} if $E^1 (e)$ has rank $m$. \smallskip  

We begin our study with the case of real spin factors.

\begin{theorem}\label{thm reduced real spin factor}
Let $W= X_{1} \mathop{\oplus}\limits^{\ell_{1}} X_{2}$ be a reduced real spin factor of type $IV_{n}^{r,s}$ with $n = r+s \geq 3$ and $ r,s \geq 1$. Suppose $I$ is a weak$^*$-closed subtriple of $W$. Then the following are equivalent:
    \begin{enumerate}[$(a)$]
    \item $I$ is an inner ideal.
    \item $I$ is weak$^*$-Hahn--Banach smooth.
\end{enumerate}
\end{theorem}

\begin{proof} $(a)\Rightarrow (b)$ follows from Corollary~\ref{c characterization of weak* inner ideals}.\smallskip
	
$(b)\Rightarrow (a)$ It is well-known that $W$ has rank 2. Fix a non-zero tripotent $e \in I$. We have to possibilities: $r(e) =2,$ and thus $e$ is a complete and unitary tripotent in $W,$ or $r(e) = 1,$ and thus $e$ is a minimal tripotent in $W$.\smallskip

Let us first assume $r(e) = 2$. In this case $W = W_{2}(e)$ and $e = (e_{1}, 0)$ for some norm-one element $e_{1} \in X_{1}$ or $e = (0, e_{2})$ for some norm-one element $e_{2} \in X_{2}$. If $e = (e_{1}, 0)$ for some norm-one $e_{1} \in X_{1}$, for each $x \in X_{2}$ write $\hat{x} :=(0,x)$. Since
$$\{e,\hat{x},e\} = 2 \langle e|\hat{x}\rangle e - \langle e|\bar{e}\rangle \bar{\hat{x}} = \hat{x}.$$
So, $(0,X_{2}) \subseteq W^{1}(e) \subseteq I$. Notice that in this case $X_{2}$ is non-zero, which implies the existence of a non-zero norm-one element $e_{2} \in X_{2}$ such that $ (0,e_{2}) \subseteq I$. Therefore $f:= (0, e_{2}) \in I$ is a complete and unitary tripotent in $W$. Similar arguments to those given above assure that $(X_{1},0) \subseteq W^{1}(f) \subseteq I$. Thus $W_{2}(e) =W = X_{1} \mathop{\oplus}\limits^{\ell_{1}} X_{2} = I$. The case $e = (0, e_{2})$ for some norm-one element $e_{2} \in X_{2}$ follows similarly. Consequently, if $I$ contains a complete tripotent, then $I=W$.\smallskip 

An alternatively proof can be given as follows. Let us write $I = I_1 \mathop{\oplus}\limits^{\ell_{1}} I_{2}$, where each $I_j$ is a closed subspace of $X_j$. Assume as before that $I$ contains an element of the form $(e_1,0)$ for a norm-one element $e_1\in X_1$. We claim that $\{0\} \mathop{\oplus}\limits^{\ell_{1}} X_{2}\subseteq I$ (equivalently, $I_2 = X_2$). Otherwise, there exists  $z_2\in X_2$ such that $\|z_2\|= 1$ and $\langle z_2| I_2\rangle =\{0\}$. For $t\in [-1,1]$, define $\varphi\in I_*$ and $\phi_t\in W_* = X_1 \mathop{\oplus}\limits^{\ell_{\infty}} X_{2}$ given by $\varphi (x) =  \langle x | e_1\rangle$ $(x\in I)$, and $\phi_t (x) =  \langle x | e_1+ t z_2\rangle$ $(x\in W),$ respectively. By construction, $\|\varphi\| = \|\phi_t\| = 1$ for all $t \in [-1,1]$ and $\phi_t |_{I} = \varphi$, which contradicts Corollary~\ref{c to Phelps theorem for weakstar topology}. This proves the claim and the rest follows as above.\smallskip

We may finally assume, without loss of generality, that every tripotent in $I$ is a minimal tripotent. Since $W$ is reduced, we can conclude that $W^{-1}(e) = \{0\}$ for all non-zero tripotent $e \in I$ (cf. Remark~\ref{remark minimal tripotents in Hahn--Banach smooth subtriples}). Corollary~\ref{c characterization of weak* inner ideals} implies that $I$ is an inner ideal of $W$.
\end{proof}

According to \cite{FriedmanRusso1985structure} and \cite{PeStacho2001}, a complex or real JBW$^*$-triple $W$ is called \emph{atomic}\label{def atomic} if it coincides with the weak$^*$-closure of the linear span of its minimal tripotents. Let $\mathcal{W}$ denote the complexification of $W$, and let $\tau$ be a conjugation (in particular a weak$^*$-continuous conjugate-linear triple automorphism of period-$2$) on $\mathcal{W}$ satisfying $\mathcal{W}^{\tau} = W$. It is almost explicit in \cite{PeStacho2001} that $W$ is atomic if, and only if, $\mathcal{W}$ is atomic. Actually, $\mathcal{W}$ is atomic if, and only if, it writes the $\ell_{\infty}$-sum of a family of Cartan factors $\{C_i: i\in \Gamma\}$. We can write $\Gamma$ as the disjoint union of the sets $\Gamma_0=\{i\in \Gamma : \tau (C_i) = C_i\}$, $\Gamma_1$, and $\Gamma_{-1}$ such that for each $i\in \Gamma_{1}$, $\tau (C_i) = C_{\sigma(i)} \perp C_i$ for a unique $\sigma (i) \in \Gamma_{-1}$. Furthermore, $$W = \mathcal{W}^{\tau} = \left(\bigoplus_{i\in \Gamma_0}^{\infty} C_i^{\tau}\right) \bigoplus^{\infty} \left(\bigoplus_{i\in \Gamma_1}^{\infty} (C_i)_{\mathbb{R}}\right).$$ The real Cartan factors $C_i^{\tau}$ as well as the realifications of the complex Cartan factors  $(C_i)_{\mathbb{R}}$ are atomic real JBW$^*$-triples (see, for example, \cite{Kaup1997}). Reciprocally, if $W$ is atomic, having in mind that every minimal tripotent in $W$ is minimal in $\mathcal{W}$ or can be written as the orthogonal sum of two orthogonal minimal tripotents in $\mathcal{W}$ (see \cite[Lemma 3.2]{PeStacho2001}), it can be easily seen that $\mathcal{W}$ is atomic. This is an argument to show that every atomic real JBW$^*$-triple coincides with the direct sum of an appropriate family of real Cartan factors or realifications of complex Cartan factors (cf. \cite[Proof of Proposition 3.1]{PoloPe2004surjective} and \cite{PeStacho2001}).  All these arguments combined with \cite[Proposition 5.1]{LiLiuPe2024characterization} allow us to obtain the next result whose proof is left to the reader. 

\begin{proposition}\label{prop equiv of atomic}
	Let $W$ be a real JBW$^*$-triple. Then $W$ is atomic if, and only if, for every non-zero tripotent $u \in W$, there exists a minimal tripotent $v \in W$ such that $v \leq u$.
\end{proposition}

By relying on the previous result it is easy to show, via Zorn's lemma, that every non-zero tripotent $e$ in an atomic JBW$^*$-triple $W$ can be written in the form $e = \sum_{i} e_i$, where $(e_i)_i$ is a family of mutually orthogonal minimal tripotents in $W$ and the family is summable with respect to the weak$^*$-topology.

\begin{lemma}\label{lem I is atomic} Let $I$ be a (non-zero) weak$^*$-closed subtriple of an atomic real JBW$^*$-triple $W$. Suppose $I$ is weak$^*$-Hahn--Banach smooth. Then $I$ is an atomic real JBW$^*$-triple.
\end{lemma}

\begin{proof} Theorem~\ref{t characterization of closed subtriples with weak*HB smoothness as hereditary subtriples} assures that $I$ is a hereditary subtriple of $W$. Take a non-zero tripotent $u \in I$. Since $W$ is atomic, we can find a minimal tripotent $e\in W$ with $e\leq u$ (cf. Proposition~\ref{prop equiv of atomic}). The tripotent $e$ lies in $I$ since the latter is hereditary. Clearly $e$ is minimal in $I$.  A new application of Proposition~\ref{prop equiv of atomic} shows that $I$ is atomic. 
\end{proof}

Let us recall some patterns followed by Peirce subspaces associated with orthogonal tripotents in a real JB$^*$-triple $E$. Let $e_1$ and $e_2$ be two orthogonal tripotents in $E$. Since $e_i \in E_0 (e_j)$ for all $i\neq j$ in $\{1,2\}$, the Peirce projections $P_k (e_1)$ and $P_l (e_2)$ commute for all $k,l\in \{0,1,2\}$ (cf. \cite[$(1.9)$ and $(1.10)$]{Horn1987characterization} whose proof is valid inthe real setting too). It is part of the basic theory of real and complex JB$^*$-triples that the following identities hold:
\begin{equation}\label{eq Peirce2 for two orthogonal}
E_2 (e_1 + e_2)  = E_2 (e_1 ) \oplus E_2 ( e_2) \oplus \left(E_1 (e_1 )\cap E_1 ( e_2) \right),
\end{equation}
\begin{equation}\label{eq Peirce1 for two orthogonal}
	E_1 (e_1 + e_2) = \left(E_1 (e_1 )\cap E_0 ( e_2) \right) \oplus \left(E_0 (e_1 )\cap E_1 ( e_2) \right).
\end{equation} Suppose now that $e_1,\ldots, e_n$ are mutually orthogonal tripotents in $E$, and set $e = e_1+\ldots+ e_n$. By applying \eqref{eq Peirce2 for two orthogonal} and \eqref{eq Peirce1 for two orthogonal}, a simple induction argument leads to 
 \begin{equation}\label{eq Peirce2 for n orthogonal}
 	E_2 (e) = \left(\bigoplus_{i=1}^n E_2 (e_i )\right) \bigoplus \left(\bigoplus_{i\neq j\in\{1,\ldots,n\}} E_1 (e_i )\cap E_1 ( e_j) \right) \ (\forall n\in \mathbb{N}).
 \end{equation}

It is the moment to present our main conclusion concerning weak$^*$-Hahn--Banach smooth weak$^*$-closed subtriples of reduced and atomic real JBW$^*$-triples.

\begin{theorem}\label{thm reduced  atomic JBW-star triples}
    Let $W$ be a reduced and atomic real JBW$^*$-triple. Suppose $I$ is a weak$^*$-closed subtriple of $W$. Then the following are equivalent:
    \begin{enumerate}[$(a)$]
    \item $I$ is an inner ideal.
    \item $I$ is weak$^*$-Hahn--Banach smooth.
\end{enumerate}
\end{theorem}

\begin{proof} Thanks to \Cref{c characterization of weak* inner ideals}, it suffices to prove that $(b)$ implies $(a)$. Suppose $I$ is weak$^*$-Hahn--Banach smooth. Lemma~\ref{lem I is atomic} assures that $I$ is reduced and atomic. Furthermore,  by Theorem~\ref{t characterization of wstarHahnBanach smooth subtriples}, $\displaystyle I= \bigcup_{e\in \mathcal{U}(I)} W^1(e),$ and $I^{1}(e) = W^{1}(e)$ for all $e\in \mathcal{U}(I)$. Moreover, $\mathcal{U}_{min}(I) \subseteq \mathcal{U}_{min}(W)$ (cf. Remark~\ref{remark minimal tripotents in Hahn--Banach smooth subtriples}).\smallskip
	
If $e$ is a minimal tripotent in $I$, we know from the fact that $W$ is reduced that $W^{-1}(e) = \{0\}\subseteq I$. Consequently, if $I$ has rank-one, every tripotent in $I$ is minimal and $W^{-1}(e) = \{0\}\subseteq I$ for every $e \in \mathcal{U}(I)$.  \Cref{c characterization of weak* inner ideals} implies that $I$ is an inner ideal.\smallskip	
	
We assume now that $I$ has rank greater than or equal to 2. Let us fix a tripotent $e \in \mathcal{U}(I)$. We can clearly assume that $e$ is not minimal by the arguments in the previous paragraph. Since $I$ is atomic, there exists a (possibly finite) family $\{u_{i}\}_{i\in \Gamma}$ of pairwise orthogonal minimal tripotents in $I$ such that $\displaystyle e = \sum_{i\in \Gamma } u_{i},$ where the series converges with respect to the weak$^*$-topology.\smallskip

\noindent \emph{Claim:} for each finite subset $\mathcal{F}\subset \Gamma$ we have $\displaystyle W_{2} \left(e_{_\mathcal{F}}\right)\subseteq I,$ where $\displaystyle e_{_\mathcal{F}}= \sum_{i\in \mathcal{F}} u_i$.\smallskip

To prove the claim, observe that by \eqref{eq Peirce2 for n orthogonal} we have 
    \begin{align}\label{thm equa W2e is the weak-star-closure of W2ui and W11}
    	W_2 (e_{_\mathcal{F}}) = \left(\bigoplus_{i=1}^n W_2 (u_i )\right) \bigoplus \left(\bigoplus_{i\neq j\in\{1,\ldots,n\}} W_1 (u_i )\cap W_1 (u_j) \right).
    \end{align} It follows from the hypotheses that $W_{2}(u_{i}) = W^{1}(u_{i}) = \mathbb{R} u_i \subset I$. We only need to prove that $W_{1}(u_{i}) \cap W_{1}(u_{j})$ is in $I$ for each $i \neq j.$ Having in mind \eqref{eq Peirce2 for two orthogonal}, it suffices to show that $W_{2}(u_{i} + u_{j})\subseteq I$ for any $i,j\in \Gamma$ with $i \neq j$.\smallskip 
    
For $i\neq j$ in $\Gamma$, Corollary 2.8 in \cite{PoloPe2004surjective} asserts that  $W_{2}(u_{i} +u_{j})$ is a real spin factor, and it is reduced since $W$ is reduced. Consider the subtriple $I \cap W_{2}(u_{i} +u_{j}) = I_{2}(u_{i} +u_{j}) \subseteq I$ which is an inner ideal of $I$. Fix any  linear functional $\phi \in I_{2}(u_{i} +u_{j})_{*}$. We may assume, without loss of generality, that $\| \phi \| =1 = \phi(s)$, where $s \in \mathcal{U}(I_{2}(u_{i} +u_{j}))$ is the support tripotent of $\phi$ in $I_{2}(u_{i} +u_{j})$. Since $I_{2}(u_{i} +u_{j})$ is an inner ideal of $I$, \Cref{c characterization of weak* inner ideals} guarantees the existence of a unique linear functional $\tilde{\phi} \in I_{*}$ such that 
    $$\tilde{\phi}\arrowvert_{I_{2}(u_{i} +u_{j})} = \phi \text{ and } \| \tilde{\phi} \| = \|\phi\| = 1.$$
Similarly, since $I$ is weak$^*$-Hahn--Banach smooth in $W$, there exists a unique linear functional $\psi \in E_{*}$ such that 
    $$\psi\arrowvert_{I} = \tilde{\phi} \text{ and } \|\psi\| =  \| \tilde{\phi} \| = 1.$$ Set $\tilde{\psi}:= \psi\arrowvert_{W_{2}(u_{i}+u_{j})} \in W_{2}(u_{i}+u_{j})_{*}$. Clearly $\tilde{\psi}$ is a  norm-preserving extension of $\phi$ in $W_{2}(u_{i}+u_{j})_{*}$. Namely,     $$\tilde{\psi}\arrowvert_{I_{2}(u_{i}+u_{j})} = \psi\arrowvert_{I_{2}(u_{i}+u_{j})} = \tilde{\phi}\arrowvert_{I_{2}(u_{i}+u_{j})} = \phi.$$ and 
    $$1 = \tilde{\psi}(s) \leq \| \tilde{\psi} \| \leq \|\psi\| = 1 = \|\phi\|.$$ We shall next show that $\tilde{\psi}$ is the unique extension of $\phi$. Suppose that $\varphi \in W_{2}(u_{i}+u_{j})_{*}$ is another norm-preserving extension of $\phi$. Since $\tilde{\psi}(s) = \|\tilde{\psi}\| =1$ and $\varphi(s) = \phi(s) = \|\varphi\| =1,$ it follows from \cite[Lemma 2.7]{PeStacho2001} that 
    $$\tilde{\psi} = \tilde{\psi} \circ P^{1}(s) \text{ and } \varphi = \varphi \circ P^{1}(s),$$
    where $P^{1}(s)$ is the corresponding Peirce projection from $W_{2}(u_{i}+u_{j})$ onto $(W_{2}(u_{i}+u_{j}))^{1}(s)$. Note that $s \in \mathcal{U}(I_{2}(u_{i} +u_{j})) \subseteq \mathcal{U}(I)$, and thus $W^{1}(s) \subseteq I$ and $(W_{2}(u_{i}+u_{j}))^{1}(s) \subseteq I_{2}(u_{i}+u_{j})$. Therefore, for every $x \in W_{2}(u_{i}+u_{j})$, we have
    \begin{align*}
        \varphi(x) &= \varphi (P^{1}(s)(x)) = \phi (P^{1}(s)(x)) = \tilde{\psi}(P^{1}(s)(x)) = \tilde{\psi}(x),
    \end{align*}
    which implies that $\varphi = \tilde{\psi}$. \smallskip

The arbitrariness of the functional $\phi \in I_{2}(u_{i} +u_{j})_*$ in the above arguments leads to the conclusion that $I_{2}(u_{i} +u_{j})$ is weak$^*$-Hahn--Banach smooth in the real (reduced) spin factor $W_{2}(u_{i} +u_{j})$. \Cref{thm reduced real spin factor} proves that $I_{2}(u_{i} +u_{j})$ is an inner ideal in $W_{2}(u_{i} +u_{j}),$ and since $u_{i} +u_{j}\in I_{2}(u_{i} +u_{j})$ we derive that $W_{2}(u_{i} +u_{j}) = I_{2}(u_{i} +u_{j})$, which finishes the proof of the Claim.\smallskip

Fix an arbitrary $x\in W_2(e)$. Let $\mathcal{F}(\Gamma)$ denote the collection of all finite subsets of $\Gamma$. We consider the net $( P_2(e_{_\mathcal{F}} ) (x) )_{_{\mathcal{F}\in \mathcal{F}(\Gamma)}}$. The previous Claim assures that $P_2(e_{_\mathcal{F}} ) (x) \in I$ for all $\mathcal{F}\in \mathcal{F}(\Gamma)$. Since the net $(e_{_\mathcal{F}} )_{_{\mathcal{F}\in \mathcal{F}(\Gamma)}}$ converges to $e$ in the strong$^*$-topology of $W$, and  $$P_2(e_{_\mathcal{F}} ) (x) = \{e_{_\mathcal{F}}, \{e_{_\mathcal{F}}, x,e_{_\mathcal{F}}\},e_{_\mathcal{F}}\},$$ it follows from the joint strong$^*$-continuity of the triple product of $W$ that $( P_2(e_{_\mathcal{F}} ) (x) )_{_{\mathcal{F}\in \mathcal{F}(\Gamma)}}\longrightarrow P_2(e) (x) =x$ in the strong$^*$-topology of $W$. Since strong$^*$-convergence implies weak$^*$-convergence, we obtain that $x\in \overline{I}^{w^*} =I$. We have proved that $W_{2}(e) \subseteq I$ for all tripotent $e\in \mathcal{U}(I)$,  which implies that $I$ is an inner ideal (cf. \Cref{c characterization of weak* inner ideals}).
\end{proof}

As in Example~\ref{example w*HahnBanach smoothness is not enough} we can find rank-one weak$^*$-closed subtriples in non-reduced real Cartan factors which are weak$^*$-Hahn-Banach smooth but not inner ideals.

\begin{example}\label{example non-reduced real Cartan factors} Suppose $W$ is an atomic non-reduced real JBW$^*$-triple. By assumption, there exists a non-zero minimal tripotent $e \in W$ satisfying $W^{1}(e) = \mathbb{R}e$ and $W^{-1}(e) \neq \{0\}$. Let $I = W^{1}(e)$ and $0\neq x \in W^{-1}(e)$. Clearly $I$ is a weak$^*$-closed subtriple of $W$ which is (weak$^*$-)Hahn-Banach smooth. However, $I$ is not an inner ideal. Observe that $I$ is a rank-one JBW$^*$-triple. 
\end{example}

Within the list of non-reduced real Cartan factors, the factors $IV_{n}^{n,0}$ and  $V^{\mathbb{O}^{rd}}$ have rank-one. So, they can be excluded if we attend to prove that every weak$^*$-closed and weak$^*$-Hahn-Banach smooth real subtriple having rank$\geq2$ is an inner ideal. This will be explored in the next sections, where we shall try to prove that the key property is the rank of the subtriple.

\section{Real weak$^*$-Hahn--Banach smooth subtriples of complex Cartan factors}\label{sec: complex Cartan factors}

We have seen in Example~\ref{example w*HahnBanach smoothness is not enough} that every complex Cartan factor admits a weak$^*$-closed real subtriple which is (weak$^*$-)Hahn-Banach smooth and is not an inner ideal. Note that all the examples in this line  have rank-one. We shall see next that no counter-example of a weak$^*$-closed real subtriple of rank bigger than or equal to $2$ can be obtained.\smallskip

We recall that a couple of tripotents $u,v$ in a (real or complex) JB$^*$-triple $E$ are said to be \emph{collinear} (written $u\top v$) if $u\in E_1(v)$ and $v\in E_1(u)$. An ordered quadruple $(e_{1},e_{2},e_{3},e_{4})$ of tripotents in a (real or complex) JB$^*$-triple $E$ is called a \emph{quadrangle} if $e_{1}\perp e_{3}$, $e_{2}\perp u_{4}$, $u_{1}\top u_{2},$ $ u_{2}\top u_{3},$ $u_{3}\top u_{4}$ $ u_{4}\top u_{1}$ and $u_{4}=2 \{{u_{1}},{u_{2}},{u_{3}}\}$ (it follows from the Jordan identity that the last equality holds if the indices are permuted
cyclically, e.g. $u_{2} = 2 \{{u_{3}},{u_{4}},{u_{1}}\}$).\smallskip

Let $u$ and $v$ be tripotents in a (real or complex) JB$^*$-triple $E$. We say that $u$ \emph{governs} $v$ ($u \vdash v$ in short) whenever $v\in U_{2} (u)$ and $u\in U_{1} (v)$. An ordered triplet $ (v,u,\tilde v)$ of tripotents in $E$, is called a \emph{trangle} if $v\bot \tilde v$, $u\vdash v$, $u\vdash \tilde v$ and $ v = Q(u)\tilde v$. 

\begin{theorem}\label{t real subtriples of complex Cartan factors} Let $I$ be a weak$^*$-closed real subtriple of a complex Cartan factor $C$. Suppose $I$ satisfies the following hypotheses:\begin{enumerate}[$(a)$]
\item $I$ is weak$^*$-Hahn-Banach smooth.
\item $I$ has rank bigger than or equal to $2$. 
\end{enumerate} Then $I$ is a complex subtriple and an inner ideal of $C$.
\end{theorem}

\begin{proof} We first observe that Lemma~\ref{lem I is atomic} assures that $I$ is an atomic real JBW$^*$-triple. According to \cite[Proof of Proposition 3.1]{PoloPe2004surjective} and \cite{PeStacho2001} (see the discussion before Proposition~\ref{prop equiv of atomic}), there is a family of non-trivial real Cartan factors or realifications of complex Cartan factors  $\{C_i\}_{i\in \Gamma}$ such that $I = \bigoplus_{i\in \Gamma} C_i$.\smallskip 
	
We claim that $\Gamma$ reduces to a single element, and hence $I$ is a real Cartan factor or a realification of a complex one. To see the claim,  let us take an arbitrary minimal tripotent $e\in I$. Clearly, there exists a unique $i_0\in \Gamma$ such that $e\in C_{i_0}$.  If there exists $i_1\in \Gamma$ with $i_1\neq i_0$, we can pick another minimal tripotent $v\in C_{i_1}\subseteq I$ which must be clearly orthogonal to $e$. By construction, $I^{1} (e+v) =\{0\}$. Remark~\ref{remark minimal tripotents in Hahn--Banach smooth subtriples} proves that $e$ and $v$ are minimal tripotents in the complex Cartan factor $C$, which for sure has rank $\geq 2$. An application of \cite[Lemma 3.10]{FerPe2018Adv} assures the existence of a non-zero tripotent $w\in C^{1} (e+v)$. Since $e+v$ is a tripotent in $I$, and the latter is a weak$^*$-Hahn-Banach smooth in $C$, Theorem~\ref{t characterization of wstarHahnBanach smooth subtriples} assures that $C^{1} (e+v)\subseteq I,$ and thus $0\neq w \in I^{1} (e+v)$, which is impossible. This finishes the proof of the claim.\smallskip
	
We have therefore shown that $I$ is a real or complex Cartan factor with rank $\geq 2$. Proposition 5.8 in \cite{Kaup1997}$(i)$ and $(ii)$ now assures the existence of a minimal tripotent $v\in I$ orthogonal to $e$. Remark~\ref{remark minimal tripotents in Hahn--Banach smooth subtriples} implies that $e$ and $v$ are minimal in $C$. An application of Lemma 3.10 \cite{FerPe2018Adv} assures that one of the next statements holds:
\begin{enumerate}[$(a)$]\item There exist minimal tripotents $v_2,v_3,v_4$ in $C$, and a unitary complex number $\delta$ such that $(e,v_2,v_3,v_4)$ is a quadrangle and $v = \delta v_3$;
\item There exist a minimal tripotent $v_3\in C$, a rank two tripotent $u\in C$, and a unitary complex number $\delta$ such that $(e, u, v_3)$ is a trangle, and $v = \delta v_3$.
\end{enumerate}

We treat both cases at once. If in case $(a)$ we set $u = v_2+v_4$, the triplet $(e, u, v_3)$ is a trangle fulfilling the conclusions in $(b)$. So, we can reduce to case $(b)$. The tripotent $e+ v = e+ \delta v_3$ lies in $I$. Take a unitary complex number $\mu$ with $\mu^2 = \delta$. It is not hard to check that the element $\mu u $ belongs to $C^{1} (e+ v)$ and the latter is contained in $I$ by Theorem~\ref{t characterization of wstarHahnBanach smooth subtriples}. Therefore, $\mu u \in I.$ Observe now that elements of the form $\alpha e + \delta \overline{\alpha} v_3$ belong to $C^{1} (\mu u )$ for all $\alpha\in \mathbb{C}$. A new application of  Theorem~\ref{t characterization of wstarHahnBanach smooth subtriples} leads to the conclusion that $\alpha e + \delta \overline{\alpha} v_3 = \alpha e + \overline{\alpha} v$ belongs to $I$ for all $\alpha\in \mathbb{C}$.  Since $e\in I$ and $I$ is a subtriple of $C$, it follows that $\alpha e = P_2 (e) (\alpha e + \delta \overline{\alpha} v_3) \in I$ for all complex number $\alpha$.\smallskip

Since in the above arguments, $e$ is an arbitrary minimal tripotent in $I$, we conclude that $\mathbb{C} e\subseteq I$ for all minimal tripotent $e\in I$. Any other tripotent $\tilde{e}\in I$ writes in the form $\tilde{e} = w^*\hbox{-}\sum_{i\in \Upsilon} u_i$, where $(u_i )_{i\in \Upsilon}$ is a family of mutually orthogonal minimal tripotents in $I$. Having in mind that $I$ is weak$^*$-closed we get from the previous conclusion that $\mathbb{C} \tilde{e}\in I$, for every tripotent $\tilde{e}\in I$.\smallskip

It is known that every element in a real JBW$^*$-triple can be approximated in norm by an algebraic element, that is, an element which can be expressed as finite linear (real) linear combination of mutually orthogonal tripotents (see the proof of $(i)\Rightarrow (ii)$ in \cite[Theorem 4.8]{IsidroKaupPalacios1995realform}),  since complex multiples of tripotents in $I$ are also in $I$, we conclude that $I$ is a complex subspace. Moreover, for each tripotent $\tilde{e}\in I$, we have $C^{-1} (e) = i C^{1} (e) \subseteq i I \subseteq I$, and thus Corollary~\ref{c characterization of weak* inner ideals} asserts that $I$ is an inner ideal of $C$.      
\end{proof}

Two non-reduced real Cartan factors with rank $\geq 2$ have been left unexplored in the previous section~\ref{sec: reduced Cartan factors}. The concrete factors are those of type $I_{2p,2q}^{\mathbb{H}}$ and $III_{2p}^{\mathbb{H}}$ ($q\geq p \geq 2$); and they both have rank $p$.%\smallskip

\begin{proposition}\label{p type I2p2q{{H}}} Let $P$ and $Q$ be Hilbert spaces of dimensions $p,q$ over the quaternion field $\mathbb{H}$ with $q\geq p \geq 2$, and let $C=\mathcal{L}(P,Q)$ be a type $I_{2p,2q}^{\mathbb{H}}$ real Cartan factor. Let $I$ be a weak$^*$-closed subtriple of $I$ which is weak$^*$-Hahn-Banach smooth in $C$ and has rank $\geq 2$. Then $I$ is an inner ideal of $C$.
\end{proposition}

\begin{proof}  Every minimal tripotent in $C$ is of the form $e= \xi\otimes \eta$, where $\xi\in Q$ and $\eta\in P$ are norm-one elements and $\xi\otimes \eta$ is the operator given by $\xi\otimes \eta (\zeta) = \langle\zeta|\eta\rangle_{_{\mathbb{H}}} \xi$ ($\zeta\in P$). Here we write $\langle \cdot | \cdot\rangle_{_{\mathbb{H}}}$ denotes the inner product in both spaces $P$ and $Q$. Let $\overline{\cdot}$ denote the conjugation on $\mathbb{H}$ defined by $\overline{\alpha_1 + \alpha_2 \vec{\imath} + \alpha_3 \vec{\jmath} + \alpha_4 \vec{k}} = \alpha_1 - \alpha_2 \vec{\imath} - \alpha_3 \vec{\jmath} - \alpha_4 \vec{k}$, for all $h = \alpha_1 + \alpha_2 \vec{\imath} + \alpha_3 \vec{\jmath} + \alpha_4 \vec{k}\in \mathbb{H}$. In this case $C^{1} (e) = \mathbb{R} e$ and $C^{-1} (e) = \{h \cdot e : h\in \mathbb{H}, h = \alpha_2 \vec{\imath} + \alpha_3 \vec{\jmath} + \alpha_4 \vec{k}\}$. \smallskip
	
As in the proof of Theorem~\ref{t real subtriples of complex Cartan factors}, $I$ is an atomic real JBW$^*$-triple (cf. Lemma~\ref{lem I is atomic}), and actually a real Cartan factor. Given a minimal tripotent $e\in I,$ we can find another minimal tripotent $v\in I$ such that $e\perp v$. Both tripotents, $e$ and $v$ are minimal in $C$ (see Remark~\ref{remark minimal tripotents in Hahn--Banach smooth subtriples}). So, we can find norm-one elements $\xi_1,\xi_2\in Q$ and $\eta_1,\eta_2\in P$ such that $e =\xi_1\otimes \eta_1$ and $v = \xi_2\otimes \eta_2$, $\xi_1\perp \xi_2$, and $\eta_1\perp \eta_2$.  The tripotent $w = \xi_2\otimes \eta_1 + \xi_1\otimes \eta_2$ belongs to $C^{1} (e+v) = I^{1} (e+v)\subseteq I$ (cf. Theorem~\ref{t characterization of wstarHahnBanach smooth subtriples}). It is easy to check that elements of the form $h \cdot e + \overline{h} \cdot v$ belong to $C^{1} (w) = I^{1} (w)\subseteq I$ for all $h\in \mathbb{H}$.  Since $I$ is a subtriple of $C$, the element $h \cdot e  = P_2 (e) 
(h \cdot e + \overline{h} \cdot v)\in I$, for all $h\in \mathbb{H}$, that is, $C_2 (e) = \mathbb{H}\cdot e\subseteq I$, for all minimal tripotent $e\in I$. Furthermore, setting $w_1 = \xi_2\otimes \eta_1 $ and $w_2= \xi_1\otimes \eta_2$, we can easily see that elements of the form $h \cdot w_1 + \overline{h} \cdot w_2$ belong to $C^{1} (e+v) = I^{1} (e+v)\subseteq I,$ for all $h\in \mathbb{H}$, and consequently, $h \cdot w_1 = P_2 (w_1) (h \cdot w_1 + \overline{h} \cdot w_2) \in I$,  for all $h\in \mathbb{H}$. Similarly, $h \cdot w_2 \in I$,  for all $h\in \mathbb{H}$. This implies that $C_1 (e+v) = C_2(w_1) \oplus C_2 (w_2) = \mathbb{H} w_1 \oplus \mathbb{H} w_2\subseteq I.$\smallskip

We have therefore proved that for every minimal tripotent $e\in I$ we have $C_2 (e)= I_2 (e) \subseteq I$, and if $v$ is any other minimal tripotent in $I$ with $e\perp v$, we also have $C_2(e+v) = I_2(e+v)\subseteq I$. We can literally repeat the arguments in the final part of the proof of Theorem~\ref{thm reduced  atomic JBW-star triples} to conclude that $C_2 (e)\subseteq I$, for all tripotent $e\in I$, and thus, $I$ is an inner ideal of $C$.    
\end{proof}

\begin{proposition}\label{p III2p{H}} Let $P$ be a Hilbert space of dimensions $p$ over the quaternion field $\mathbb{H}$ with $p \geq 2$, and let $C=\{ w \in \mathcal{L}(P): w^*=w \}$ be a type $III_{2p}^{\mathbb{H}}$ real Cartan factor. Let $I$ be a weak$^*$-closed subtriple of $I$ which is weak$^*$-Hahn-Banach smooth in $C$ and has rank $\geq 2$. Then $I$ is an inner ideal of $C$.
\end{proposition}

\begin{proof} We keep the notation in the proof of the previous Proposition~\ref{p type I2p2q{{H}}}. The arguments there can be easily adapted to our setting by just observing that each minimal tripotent in $C$ is of the form $e=\vec{k}\cdot \xi \otimes \xi,$ for a unitary vector $\xi \in P$. It is not hard to see that $C^{1} (e) =\mathbb{R} \cdot e$, $C^{-1} (e) = \mathbb{H}_{1,\vec{k}}^{\perp} e$, and $C_2 (e) = \mathbb{H}_{\vec{k}}^{\perp} e$, where $\mathbb{H}_{1,\vec{k}}^{\perp} = \hbox{span}_{\mathbb{H}} \{\vec{\imath},\vec{\jmath}\}$, and $\mathbb{H}_{\vec{k}}^{\perp} = \hbox{span}_{\mathbb{H}} \{1,\vec{\imath},\vec{\jmath}\}$. For the sake of brevity, the rest is left to the reader. 
\end{proof}

\section{Applications: A sufficient condition to be an inner ideal}\label{sec: conclusions}

In order to present our main conclusion concerning closed subtriples which are Hahn-Banach smooth, we recall first the remaining aspects in the Gelfand-Naimark theorem for real JB$^*$-triples (see \cite[Proposition 3.1]{PoloPe2004surjective} and \cite{PeStacho2001}). Let $E$ be a real JB$^*$-triple. It is well-known that we can naturally embed $E$ into its bidual space $E^{**}$ and the latter is a real JBW$^*$-triple \cite{IsidroKaupPalacios1995realform}. Let $\iota_{_E} : E\hookrightarrow E^{**}$ the natural isometric triple embedding with weak$^*$-dense image. The atomic decomposition of $E^{**}$ decomposes the latter as the direct sum of two weak$^*$-closed (orthogonal) triple ideals $A$ and $N$, $A$ being the weak$^*$-closed real linear span of all the minimal tripotents in $E^{**}$ while $N$ contains no minimal tripotents \cite[Theorem 3.6]{PeStacho2001}. The quoted result (see also the proof of \cite[Proposition 3.1]{PoloPe2004surjective}) also assures that if $\pi_{_{E^{**}_{at}}}: E^{**}\to A$ denotes the natural projection of $E^{**}$ onto $A$, the mapping $\pi_{_{E^{**}_{at}}} \circ \iota_{_E}: E\hookrightarrow A$ is an isometric triple embedding with weak$^*$-dense image, and $\displaystyle A =\bigoplus^{\infty}_{i\in \Gamma_{_{E^{**}}}} C_i,$ where each $C_i$ is a real Cartan factors or a realification of a complex Cartan factor and $\Gamma_{_{E^{**}}}$ is an index set. For each ${i\in \Gamma_{_{E^{**}}}}$ the natural projection of $E^{**}$ onto $C_i$ will be denoted by $\pi^{\tiny{E}}_i$ and will be called the projection of $E^{**}$ onto its $C_i$-component.\smallskip

The atomic decomposition can be actually established for every real or complex JBW$^*$-triple (cf. \cite{FriedmanRusso1985structure,PeStacho2001}). According to what we survey in page~\pageref{def atomic}, a real (resp., complex) JBW$^*$-triple is atomic if it coincides with its atomic part, in such a case, it can be written as a direct sum of real Cartan factors or realifications or complex Cartan factors (resp., complex Cartan factors).\smallskip  

After setting the basic terminology, we can now state a result providing sufficient conditions to guarantee that a Hahn-Banach smooth closed subtriple of a real JB$^*$-triple is an  inner ideal. 

\begin{theorem}\label{t sufficient conditions for HBS implies inner ideal} Let $I$ be a closed subtriple of a real JB$^*$-triple $E$. Suppose that $I$ satisfies the following hypotheses:
\begin{enumerate}[$(a)$]\item $I^*$ is separable.
\item $I$ is Hahn-Banach smooth.
\item The projection of $I^{**}$ onto each real or complex Cartan factor summand in the atomic part of $E^{**}$ is zero or has rank greater than or equal to $2$.
\end{enumerate} Then $I$ is an inner ideal of $E$. 
\end{theorem}  

\begin{proof} We keep the notation at the beginning of this section. Observe that $I^{**} = \overline{I}^{w^*}$ is weak$^*$-Hahn-Banach smooth in $E^{**}$, which subsequently implies that $\pi^{\tiny{E}}_i (I^{**})$ is a weak$^*$-closed weak$^*$-Hahn-Banach smooth subtriple of $C_i$ ($\forall i \in \Gamma_{_{E^{**}}}$).\smallskip
	
By hypothesis, for each $i \in \Gamma_{_{E^{**}}}$, the weak$^*$-subtriple $\pi^{\tiny{E}}_i (I^{**})$ is zero or has rank $\geq 2$. We can clearly assume that it is non-zero. Having in mind that $\pi^{\tiny{E}}_i (I^{**})$ is a weak$^*$-Hahn-Banach smooth subtriple of $C_i$, and the latter is one of the following factors: a reduced Cartan factor, the realification of a complex Cartan factor, a real Cartan factor of type $I_{2p,2q}^{\mathbb{H}}$, or a real Cartan factor of type $III_{2p}^{\mathbb{H}}$ ($q\geq p\geq 2$), Theorems~\ref{thm reduced  atomic JBW-star triples} and \ref{t real subtriples of complex Cartan factors}, and Propositions~\ref{p type I2p2q{{H}}} and \ref{p III2p{H}} prove that $\pi^{\tiny{E}}_i (I^{**})$ is an inner ideal of $C_i$ for all $i \in \Gamma_{_{E^{**}}}$. It is clear, from definitions, that $\displaystyle \pi_{_{E^{**}_{at}}}(I^{**}) = \bigoplus^{\infty}_{i \in \Gamma_{_{E^{**}}}} \pi^{\tiny{E}}_i (I^{**})$, which now gives that $\pi_{_{E^{**}_{at}}}(I^{**})$ %$= \overline{\pi_{_{E^{**}_{at}}}(\iota_{_E} (I))}^{w^*}$
 is an inner ideal in the atomic part $\displaystyle A =\bigoplus^{\infty}_{i\in \Gamma_{_{E^{**}}}} C_i,$ of $E^{**}$.\smallskip

Consider now the isometric triple embedding $(\pi_{_{E^{**}_{at}}} \circ \iota_{_E}) : E\hookrightarrow A$. Given $a\in I$, and $b\in E$, the triple product $$(\pi_{_{E^{**}_{at}}} \circ \iota_{_E}) (\{a,b,a\}) = \{(\pi_{_{E^{**}_{at}}} \circ \iota_{_E}) (a), (\pi_{_{E^{**}_{at}}} \circ \iota_{_E}) (b), (\pi_{_{E^{**}_{at}}} \circ \iota_{_E}) (a) \}$$ lies in $\pi_{_{E^{**}_{at}}}(I^{**})$ because $(\pi_{_{E^{**}_{at}}} \circ \iota_{_E}) (a) \in \pi_{_{E^{**}_{at}}}(I^{**})$ and the latter is an inner ideal of $A$. We can therefore conclude that \begin{equation}\label{eq density 0609} (\pi_{_{E^{**}_{at}}} \circ \iota_{_E}) (\{a,b,a\}) \in \pi_{_{E^{**}_{at}}}(I^{**}) .
\end{equation}\smallskip

It is well known in functional analysis that $I^*$ separable is equivalent to the metrizability of the closed unit ball of $I^{**}$ when equipped with the weak$^*$ topology (cf. \cite[Proposition 3.24]{FabianHabalaMontesinosHajekBook2001}). It follows from our assumptions that $I$ is sequentially weak$^*$-dense in its bidual, and thus, by \eqref{eq density 0609}, that we can find a sequence $(a_n)_n\subseteq I$ converging in the $\sigma(I^{**},I^*)$-topology to some $\tilde{a}\in I^{**}$ such that \begin{equation}\label{eq one 0709}
(\pi_{_{E^{**}_{at}}} \circ \iota_{_E}) (\{a,b,a\}) =  \pi_{_{E^{**}_{at}}} (\tilde{a}) =\sigma(E^{**},E^{*})\hbox{-}\lim_{n} \pi_{_{E^{**}_{at}}} ( \iota_{_E} (a_n)).
\end{equation}

We need to go deeper on the geometric properties of the atomic decomposition. Let the symbol $\partial_{e} (\mathcal{B}_{E^*})$ stand for the set of all extreme points of the closed unit ball of $E^{*}$. For each $\varphi\in \partial_{e} (\mathcal{B}_{E^*})$, there exists a unique minimal tripotent $e\in E^{**}$ such that $\varphi = \varphi \circ P^1 (e)$, or equivalently, $P^{1} (e) (x) = \varphi (x) \cdot e$ for all $x\in E^{**}$ (cf. \cite[Lemma 2.7 and Corollary 2.1]{PeStacho2001}). In particular $\varphi = \varphi \circ \pi_{_{E^{**}_{at}}}$ for all $\varphi \in \partial_{e} (\mathcal{B}_{E^*})$. It follows from this fact, the conclusion in \eqref{eq one 0709}, and the basic properties of the mapping $\iota_{_E}$ that $$\varphi\{a,b,a\} = \varphi (\pi_{_{E^{**}_{at}}} \circ \iota_{_E}) (\{a,b,a\}) =  \lim_{n} \varphi \pi_{_{E^{**}_{at}}} ( \iota_{_E} (a_n)) =  \lim_{n} \varphi (a_n),$$ for all $\varphi\in \partial_{e} (\mathcal{B}_{E^*})$. An application of Rainwater's theorem asserts that the sequence $(a_n)_n$ converges weakly in $E$ to $\{a,b,a\}$. Therefore $\{a,b,a\}$ lies in the weak closure of $I$, and thus in $I$ since the latter is a closed subspace (cf. Mazur theorem \cite[Theorem 3.19]{FabianHabalaMontesinosHajekBook2001}). 
\end{proof}

The above result is new even in the case of real C$^*$- and real JB$^*$-algebras (see \cite{BNPS2025MidealsinrealJB,BNPSMidealsOperatorAlgebras} and references therein for the detailed definitions). The arguments in the proof of Theorem~\ref{t sufficient conditions for HBS implies inner ideal} can be easily adapted to get the following results.

\begin{theorem}\label{t inner ideal in atomic complex} Let $I$ be a weak$^*$-closed real subtriple of an atomic complex JBW$^*$-triple $E$. Suppose that $I$ satisfies the following hypotheses:
	\begin{enumerate}[$(a)$]
		\item $I$ is weak$^*$-Hahn-Banach smooth. 
		\item The projection of $I$ onto each Cartan factor component of $E$ is zero or has rank $\geq 2$. 
	\end{enumerate} Then $I$ is a complex JBW$^*$-subtriple and an inner ideal of $E$.
\end{theorem}

\begin{theorem}\label{t inner ideal in atomic real} Let $I$ be a weak$^*$-closed real subtriple of an atomic real JBW$^*$-triple $E$. Suppose that $I$ satisfies the following hypotheses
	\begin{enumerate}[$(a)$]
		\item $I$ is weak$^*$-Hahn-Banach smooth. 
		\item The projection of $I$ onto each Cartan factor component of $E$ is zero or has rank $\geq 2$. 
	\end{enumerate} Then $I$ is an inner ideal of $E$.
\end{theorem}

\section{Post-credits section: \texorpdfstring{$M$}{M}-ideals in real JB\texorpdfstring{$^*$}{*}-triples revisited}\label{sec: post-credits}
 
In a recent contribution, D. Blecher, M. Neal and the second and third authors of this note culminated the algebraic characterization of $M$-ideals in real JB$^*$-triples as closed subtriples (see \cite{BNPS2025MidealsinrealJB}). A different, and somehow shorter, argument valid only for real C$^*$- and JC$^*$-algebras is presented in \cite{BNPSMidealsOperatorAlgebras}. In this section we shall present a very simplified proof valid for all real JB$^*$-triple which just uses a geometric tool from \cite{BNPS2025MidealsinrealJB} and shortens several pages the arguments --we can actually avoid the results in sections 4 and 5 of the just quoted paper.

\begin{theorem}\label{t every M-summand is a weak*-closed subtriple} Let $W$ be a real JBW$^*$-triple. Then, for every  $M$-projection on $W$ the elements $P(x)$ and $ (Id-P)(y)$ are orthogonal in $W$ for all $x,y\in W$. Furthermore, every $M$-summand in $W$ is a weak$^*$-closed triple ideal of $W$.
 \end{theorem}

\begin{proof} Let $P:W\to W$ be an $M$-projection. Consider the  $M$-summand $M= P(W)$. We can clearly assume that $\{0\}\neq M\neq W$. As shown in \cite[Remark 3.1]{BNPS2025MidealsinrealJB}, it follows from  Proposition 2.3 in \cite{MartinezPe2000separate} (alternatively, from \cite[Theorem I.1.9]{HarmandWernerWernerBook1993}) that every $M$-summand in $W$ is weak$^*$-closed. Consequently, $M$ is a weak$^*$-closed subspace. We further know that there exists a weak$^*$-continuous $M$-projection $P$ of $W$ onto $M$ such that $W = M \mathop{\oplus}\limits^{\ell_\infty} N$, where $N=(Id-P) (M)$ is clearly weak$^*$-closed too.\smallskip

 The connections between $M$-summand in $W$ and the facial structure of the closed unit ball of $W$ are exploited in \cite[Proposition 3.2, Theorem 3.1 and  Corollary 3.2]{BNPS2025MidealsinrealJB}. We just need the following result in \cite[Corollary 3.2]{BNPS2025MidealsinrealJB}: for every $v\in \partial_{e}(\mathcal{B}_{M})$ and every $w\in \partial_{e} (\mathcal{B}_{N})$, the elements $v,$ $w,$ and $v \pm w$ are tripotents in $W$, $v$ is orthogonal to $w$, and $v \pm w$ is a complete tripotent in $W$.\smallskip

 We shall first prove that $M\perp N$ in $W$, that is: 
 \begin{equation}\label{eq M and N are orthogonal} \hbox{Given } m\in M, \ n\in N, \hbox{ we have } m\perp n \text{ (equivalently, } \{m,m,n\}=0 ). 
 \end{equation} Namely, since $M$ is weak$^*$-closed, it is a dual Banach space, and thus, by the Krein--Milman theorem, its closed unit ball, $\mathcal{B}_{M}$, coincides with the weak$^*$-closure of the convex hull of $\partial_e(\mathcal{B}_{M})$. The same statement holds for $\mathcal{B}_{N}$. So, given $m\in \mathcal{B}_{M}$ and $n\in \mathcal{B}_{N}$, we can approximate $m$ and $n$, in the weak$^*$-topology of $W$ by elements of the form $\displaystyle \sum_{j=1}^{k_1} t_j v_j$ and $\displaystyle \sum_{k=1}^{k_2} s_{k} w_{k}$, where $v_j\in \partial_e (\mathcal{B}_{M})$, $w_{k}\in \partial_e (\mathcal{B}_{N})$, $t_j,s_{k}\in \mathbb{R}^{+}$ for all $j,k$, $\displaystyle \sum_{j=1}^{k_1} t_j =1$, and $\displaystyle \sum_{k=1}^{k_2} s_{k} =1$. It follows from the conclusion in the second paragraph of this proof that $v_j\perp w_{k}$ for all $j,k$, and hence $$\left\{ \sum_{j=1}^{k_1} t_j v_j, \sum_{j=1}^{k_1} t_j v_j,  \sum_{k=1}^{k_2} s_{k} w_{k}\right\} =0.$$ By applying that the triple product of $W$ is separately weak$^*$-continuous (cf. \cite[Theorem 2.11]{MartinezPe2000separate}), we obtain $\{m,m,n\} =0$ for all $m\in \mathcal{B}_{M}, n\in \mathcal{B}_{N},$ and consequently for all $m\in M,$ $n\in N.$\smallskip

% Let us next show that $M$ is a JB$^*$-subtriple of $W$. Take $m_1,m_2,m_3\in M.$ By hypotheses, $\{m_1,m_2,m_3\} = m_4+ n_1$ with $m_4\in M$ and $n_1\in N.$ By the Jordan identity and \eqref{eq M and N are orthogonal} we obtain $$\begin{aligned} \{n_1,n_1,n_1\} &= \{ \{m_1,m_2,m_3\},n_1,n_1 \} = \{ \{m_1,n_1,n_1 \},m_2,m_3\} \\ &- \{m_1,\{m_2,n_1,n_1 \},m_3\} +\{m_1,m_2,\{m_3,n_1,n_1 \}\} = 0,     \end{aligned}$$ which implies that $n_1=0,$ and thus $\{m_1,m_2,m_3\} = m_4\in M.$\smallskip
 
Suppose finally that we take $m\in M$, $a,b\in W$, where we write $a = P(a) + (I-P) (a),$ and $b = P(b) + (I-P) (b)$ --with $P(a),P(b)\in M$, $(I-P) (a), (I-P)(b) \in N$. It follows from the properties proved in the above paragraph that $\{a,b,m\} = \{P(a),P(b),m\} \in M,$ and $\{a,m,b\} = \{P(a),m,P(b)\} \in M,$ which shows that $M$ is a JB$^*$-subtriple and a triple ideal. 
\end{proof}

The algebraic characterization of norm-closed $M$-ideals in a real JB$^*$-triple as obtained in \cite[Theorem 5.2]{BNPS2025MidealsinrealJB} can be easily deduced from our previous theorem.

\begin{corollary}{\rm(compare \cite[Theorem 5.2]{BNPS2025MidealsinrealJB})} Let $E$ be a real JB$^*$-triple, and let $M$ be a norm-closed subtriple of $E$. Then $M$ is an $M$-ideal of $E$ if, and only if, it is a norm-closed triple ideal of $E$.
 \end{corollary}
 
 \medskip

\textbf{Acknowledgements}\quad 
First author is partly supported by NSF of China (Grant No. 12571143). 
Second author supported by MCIN (Ministerio de Ciencia e Innovación, Spain)/AEI/ 10.13039/ 501100011033 and “ERDF A way of making Europe”  grant PID2021-122126NB-C31, Junta de Andalucía grant FQM375,  MOST (Ministry of Science and Technology of China) project in China grant G2023125007L, and by the IMAG--Mar{\'i}a de Maeztu grant CEX2020-001105-M/AEI/10.13039/501100011033. Third author supported by the China Scholarship Council (Grant No.202306740016), grant PID2021-122126NB-C31 funded by MICIU/AEI/10.13039/501100011033, and the \linebreak Shanghai Natural Science Foundation under the Science and Technology Innovation Action Plan 2024 (Grant No.24ZR1415600).\smallskip

\noindent Part of this work was completed during visits of A.M. Peralta and J. Zhang to the Chern Institute of Mathematics at Nankai University and Universidad de Granada, respectively. Both authors acknowledge their gratitude for the hospitality. The visit of to Universidad de Granada was supported by a Scholarship of Chern Class-Shiing-Shen Chern International Exchange Award.\medskip

%---------------------------------

\subsection*{Data availability}

There are no data associates for the submission entitled ``Uniqueness of Hahn--Banach extensions determines inner ideals in real C$^*$-algebras and real JB$^*$-triples''.

\subsection*{Statements and Declarations}

The authors declare that they have no financial or conflict of interest.

%\reference here
%\bibliographystyle{siam}
%\bibliographystyle{unsrt}
%\bibliographystyle{abbrv}
%\bibliography{sn-a}

%\bibliographystyle{plain}
%\bibliography{bibl062025}

% ------------------------------------------------------------------------
\end{document}